\documentclass[12pt]{amsart}
\usepackage{graphics,amssymb,latexsym}
\input epsf                 

\numberwithin{equation}{section}

\newtheorem{theorem}{Theorem}[section]
\newtheorem{proposition}[theorem]{Proposition}
\newtheorem{conjecture}[theorem]{Conjecture}
\newtheorem{corollary}[theorem]{Corollary}
\newtheorem{lemma}[theorem]{Lemma}

\theoremstyle{definition}

\newtheorem{remark}[theorem]{Remark}
\newtheorem{example}[theorem]{Example}
\newtheorem{definition}[theorem]{Definition}

\newtheorem{prop}[theorem]{Proposition}

\newcommand{\set}[1]{{\left\lbrace #1 \right\rbrace}}

\newcommand{\br}[1]{{\langle #1 \rangle}}

\newcommand{\ep}{\varepsilon}

\newcommand{\Pgem}{{\Phi_{\ge -1}^m}}

\newcommand{\Poly}{{\mathbf{P}}}
\newcommand{\overunder}[2]{
\!\begin{array}{c}
\scriptstyle{#1}\\[-.1in]
-\!\!\!-\!\!\!-\\[-.1in]
\scriptstyle{#2}
\end{array}
\!
}

\DeclareRobustCommand{\SkipTocEntry}[4]{}

\begin{document}

\title[Generalized cluster complexes 
and Coxeter combinatorics]{Generalized cluster complexes\\[.1in]
and Coxeter combinatorics}

\author{Sergey Fomin}
\address{Department of Mathematics, University of Michigan,
Ann Arbor, MI 48109, USA} 
\email{fomin@umich.edu, nreading@umich.edu}

\author{Nathan Reading}

\date{May 4, 2005.  Revised July 25, 2005}

\thanks{This work was partially supported by NSF grants DMS-0245385 (S.F.) 
and DMS-0202430 (N.R.).}\




\begin{abstract}
We introduce and study a family of simplicial complexes associated to an arbitrary
finite root system and a nonnegative integer parameter~$m$. For $m=1$, our
construction specializes to the 
(simplicial) generalized associahedra or, equivalently, to the cluster complexes for 
the cluster algebras of finite type. 

Our  computation of the face numbers and $h$-vectors of these complexes
produces the enumerative invariants defined in other contexts by
C.~A.~Athanasiadis, suggesting links to a host 
of well studied problems in algebraic combinatorics of finite Coxeter
groups, root systems, and hyperplane arrangements. 

Recurrences
satisfied by the face numbers of our complexes lead to combinatorial
algorithms for determining Coxeter-theoretic invariants. That is,
starting with a Coxeter diagram of a finite Coxeter group, one can
compute the Coxeter number, the exponents, and other classical
invariants by a recursive procedure that only uses most basic
graph-theoretic concepts applied to the input diagram. 

\end{abstract}

\maketitle


\setcounter{tocdepth}{1}
\tableofcontents

\section{Introduction}

In the first part of this paper
(Sections~\ref{sec:cluster-complexes}--\ref{sec:recip}), 
we introduce and study a simplicial complex $\Delta^m(\Phi)$ 
associated to a finite root system $\Phi$ and a 
nonnegative integer parameter~$m$. 
For $m\!=\!1$, our construction specializes to the (simplicial) generalized
associahedra  $\Delta(\Phi)$ introduced in~\cite{ga} and identified in~\cite{ca2} as
the underlying complexes for the cluster algebras of finite type.
We enumerate the faces of the complexes~$\Delta^m(\Phi)$
and determine their Euler characteristics. 
For the classical types in the Cartan-Killing classification, 
we provide explicit combinatorial descriptions of these complexes 
in terms of dissections of a convex polygon into smaller polygons. 
In types $A_n$ and~$B_n\,$, we rediscover the constructions given by
E.~Tzanaki~\cite{Tzanaki}. 
Enumerative combinatorial invariants of the complexes $\Delta^m(\Phi)$
provide natural generalizations of the Fuss-Catalan, 
Kirkman-Cayley, and Przytycki-Sikora numbers to arbitrary types,
suggesting connections with a host of well studied problems in the
algebraic combinatorics of finite Coxeter groups, root systems, and
hyperplane arrangements.

The second part of the paper (Sections~\ref{sec:coxeter-inv}--\ref{sec:fake}) 
is devoted to combinatorial algorithms
for determining Coxeter-theoretic invariants. 
Starting with a Coxeter diagram of a finite Coxeter group
(or with the corresponding Dynkin diagram or Cartan matrix),
we compute the Coxeter number, the exponents, and other related 
invariants by a procedure (in fact, by several alternative procedures) 
which is entirely combinatorial in nature. 
That is, these procedures only use the most elementary graph-theoretic
concepts applied to the input diagram,
and do not involve, in any way, the
Coxeter group itself, the associated root system, the root lattice, or
any other group-theoretic, lattice-theoretic,  or Lie-theoretic
notions and constructions.
The crucial ingredients in all of these calculations are some identities
for the face numbers of the generalized cluster complexes established
in the first part of the paper. 
Formally extending these procedures,
one can calculate ``fake'' Coxeter invariants 
(such as a ``fake Coxeter number'') of various infinite Coxeter groups. 

The only prerequisites for this paper are the fundamentals of root systems and
finite Coxeter groups; see, e.g., \cite[Lectures~1--2]{pcmi}, or the
standard references~\cite{Bourbaki, Humphreys}.


The paper is organized as follows. 
Section~\ref{sec:cluster-complexes} reviews the main notions 
related to the 
generalized associahedra~$\Delta(\Phi)$. 
Readers familiar with the original source \cite[Section~3]{ga} 
or with the survey in \cite[Section~4.3]{pcmi} 
may proceed directly to Section~\ref{sec:gcc}. 
There, we define the generalized cluster
complexes~$\Delta^m(\Phi)$ and describe their basic structural
properties. The proofs of two of these results 
(Theorems~\ref{th:alt def}~and~\ref{th:restrict})
are given in Sections~\ref{sec:alt def} and~\ref{sec:restrict},
respectively, and can be skipped at first reading. 
Section~\ref{sec:examples} presents some examples of the
complexes~$\Delta^m(\Phi)$. 
Combinatorial models for $\Delta^m(\Phi)$ when $\Phi$ is 
of classical type (i.e., $A_n\,$, $B_n\,$, $C_n\,$, or~$D_n$)
are given in Section~\ref{sec:classical}. 
Recursions and explicit multiplicative formulas for the face numbers
of 
$\Delta^m(\Phi)$ 
are given in Section~\ref{sec:face-numbers};
these formulas are proved in Section~\ref{sec:proof-product-f}. 
In Section~\ref{sec:h-vectors}, we compute the $h$-vector 
of 
$\Delta^m(\Phi)$, recovering the $m$-Narayana
numbers introduced by C.~A.~Athanasiadis~\cite{Ath-TAMS}. 
The reduced Euler characteristic of $\Delta^m(\Phi)$ is determined in
Section~\ref{sec:euler}. 
In Section~\ref{sec:recip}, we enumerate the $m$-analogues of
``positive clusters.'' \linebreak[3]
Section~\ref{sec:coxeter-inv} presents several alternative algorithms
for computing classical invariants of a root system~$\Phi$
(or the associated Coxeter group~$W$). 
In Section~\ref{sec:fake}, these algorithms are applied to calculate
``fake invariants'' of some infinite Coxeter groups.





\section{Cluster complexes of finite type} 
\label{sec:cluster-complexes}

This section reviews the basic facts from \cite[Section~3]{ga}
concerning (simplicial) generalized associahedra; alternatively, see  
\cite[Section~4.3]{pcmi}. 
The only difference is that we do \emph{not} assume that the
underlying root system~$\Phi$ is crystallographic. This does not in
fact create any additional complications: as noted in~\cite{pcmi}, 
the constructions in~\cite{ga} extend verbatim to the
non-crystallographic case. 

Let $\Phi$ be a finite root system of rank~$n$. 
We denote by $\Phi_{>0}$ the set of positive
roots in~$\Phi$. 
The set of simple roots in $\Phi$ is denoted by 
$\Pi=\{\alpha_i\,:\,i\in I\}$,
where $I$ is an $n$-element indexing set. 
(The standard choice is $I=\{1,\dots,n\}$.) 
Accordingly, $-\Pi=\{-\alpha_i\,:\,i\in I\}$ is the set of negative
simple roots. 
%
The set 
$S=\{s_i\,:\,i\in I\}$ 
of simple reflections corresponding to the simple
roots~$\alpha_i$ generates a finite reflection group~$W$. 
The pair $(W,S)$ is a Coxeter system; 
the Coxeter group~$W$ naturally acts on the set of
roots~$\Phi$. 
 
Let us temporarily assume that the root system $\Phi$
is \emph{irreducible}. 
Let $I=I_+\cup I_-$ be a decomposition of~$I$ 
such that the sets $I_+$ and $I_-$ are disjoint, 
and each of them labels a totally
disconnected set of vertices in the Coxeter diagram of~$\Phi$. 

The ground set for the cluster complex $\Delta(\Phi)$ is the set
\[
\Phi_{\geq-1} = \Phi_{>0} \cup (-\Pi) 
\]
of \emph{almost positive roots}. 
We define the involutions $\tau_\pm:\Phi_{\geq
  -1}\to\Phi_{\geq -1}$ by
\begin{equation*} 
\label{eq:tau-pm-on-roots} 
\tau_\varepsilon(\alpha) = 
\begin{cases} 
\displaystyle 
\ \ \alpha & \text{if $\alpha = - \alpha_i\,$, for $i \in I_{- \varepsilon}$;}
\\[.1in] 
\displaystyle\Bigl(\prod_{i \in I_\varepsilon} s_i\Bigr)\,(\alpha) & \text{otherwise,} 
\end{cases} 
\end{equation*} 
for $\varepsilon\in\{+,-\}$. 
The product $R=\tau_-\tau_+$ can be viewed as a deformation of
the Coxeter element in~$W$. 
We denote by $\langle R \rangle$ the cyclic group generated by~$R$. 

Let $h$ denote the \emph{Coxeter number} of~$W$ and let $w_\circ$ be the longest element of $W$.


\begin{lemma}[{\cite[Theorem~2.6, Proposition~2.5]{ga}}]
\label{lem:dihedral}
The order of $R$ is $(h+2)/2$ if \linebreak[3]
$w_\circ = -1$, and is $h+2$ otherwise.
Every $\langle R \rangle$-orbit in  $\Phi_{\geq - 1}$ has a
nonempty intersection with $(- \Pi)$.
These intersections are precisely the $\langle - w_\circ \rangle$-orbits
in~$(- \Pi)$.
\end{lemma}


The following theorem is a reformulation of results in
\cite[Section~3.1]{ga}. 

\begin{theorem} 
\label{th:compat}
There is a unique symmetric binary relation on
$\Phi_{\geq -1}$ (called ``compatibility'') such that: 
\begin{itemize}
\item
$\alpha$ and $\beta$ are compatible if and only if $R(\alpha)$
and $R(\beta)$ are compatible; 
\item
a negative simple root $-\alpha_i$ is compatible with a positive root~$\beta$ 
if and only if the simple root expansion of~$\beta$ does not
involve~$\alpha_i$. 
\end{itemize}
\end{theorem}

Following \cite[p.~983]{ga}, we define the \emph{cluster complex}
$\Delta(\Phi)$ (of type~$\Phi$) 
as the clique complex for the compatibility relation. 
That is, a subset of roots in
$\Phi_{\geq -1}$ forms a simplex in $\Delta(\Phi)$ 
if and only if every pair of roots in this subset is compatible. 

If $\Phi$ is \emph{reducible}, with irreducible
components $\Phi_1,\dots,\Phi_l$, then 
\[
\Phi_{\geq -1} = \bigcup_j (\Phi_j)_{\geq -1}
\]
(disjoint union). 
We declare two roots in $\Phi_{\geq -1}$  
compatible if and only if they either belong to different components,
or belong to the same component and are compatible within it. 
Thus, the simplicial complex $\Delta(\Phi)$ is the 
\emph{join} of the complexes $\Delta(\Phi_j)$.

The simplicial complex $\Delta(\Phi)$ is homeomorphic to a
sphere~\cite{ga}. Moreover, it can be explicitly realized as a
boundary of a convex polytope~\cite{cfz}, a polar dual to the
(simple) \emph{generalized associahedron} of type~$\Phi$. 
This is why  $\Delta(\Phi)$ is sometimes referred to as the
\emph{simplicial} generalized associahedron (of type~$\Phi$). 

By \cite[Theorem~1.13]{ca2},
the definition of $\Delta(\Phi)$ given above is equivalent to the
algebraic definition~\cite{ca2} of a cluster complex for a cluster
algebra of finite type. Although inspired by cluster algebra theory,
this paper does not rely on any of its results; cf.\
Remark~\ref{rem:m-cluster-algebras}.

\section{Generalized cluster complexes} 
\label{sec:gcc}

Let $m$ be a nonnegative integer. In this section, we define and begin to
study the main object of this paper, the generalized cluster
complex~$\Delta^m(\Phi)$. 
The ground set of~$\Delta^m(\Phi)$
is the set $\Phi^m_{\geq -1}$ of 
\emph{colored almost positive roots}. 
It consists of $m$ copies of the set $\Phi_{>0}$ of positive
roots in~$\Phi$ together with one copy of the negative simple roots.
It will be convenient to use the following notation for the elements
of $\Phi^m_{\geq -1}$. 
For  each $\alpha\in\Phi_{>0}\,$, let $\alpha^1,\dots,\alpha^m$ denote 
the $m$ ``colored'' copies of $\alpha$ occurring in~$\Phi^m_{\geq
  -1}$. 
Each negative simple root $\alpha$ occurs in~$\Phi^m_{\geq -1}$
as~$\alpha^1$. 
Thus, 
\[
\Phi^m_{\geq -1} = \{\alpha^k\,:\,\alpha\in\Phi_{>0},\, 
k\in\{1,\dots,m\}\}
\cup \{(-\alpha_i)^1\,:\,i\in I\}. 
\]

The simplicial complex~$\Delta^m(\Phi)$ is defined using 
the binary \emph{compatibility relation} on $\Phi^m_{\geq -1}$. 
This relation can be defined in (at least) two different ways, 
which we will now describe. 
As in the non-colored ($m=1$) case, we assume that $\Phi$ is
irreducible,
since the reducible case can be obtained by taking joins. 


For a root $\beta\in\Phi_{\geq -1}\,$, let $d(\beta)$ 
denote the smallest $d$ such that 
\[
\underbrace{R(R(R(\cdots R}_{d\
  {\rm times}}(\beta)\cdots )))
\]
is a negative root. In particular, $d(\beta)=0$ if $\beta$ is negative simple. 

\begin{definition}
\label{def:m-compat}
Two colored roots $\alpha^k,\beta^l\in\Phi^m_{\geq -1}$ are called
\emph{compatible} if and only if one of the following conditions is satisfied:  
\begin{itemize}
\item 
$k>l$, $d(\alpha)\leq d(\beta)$, and the roots $R(\alpha)$ and $\beta$
are compatible (in the original ``non-colored'' sense of
Theorem~\ref{th:compat});  
\item 
$k<l$, $d(\alpha)\geq d(\beta)$, and the roots $\alpha$ and $R(\beta)$
are compatible; 
\item 
$k>l$, $d(\alpha)>d(\beta)$, and  the roots $\alpha$ and $\beta$ are compatible; 
\item 
$k<l$, $d(\alpha)<d(\beta)$, and  the roots $\alpha$ and $\beta$ are compatible; 
\item 
$k=l$, and  the roots $\alpha$ and $\beta$ are compatible. 
\end{itemize}
\end{definition}

\begin{lemma}
The compatibility relation on $\Phi^m_{\geq -1}$ is symmetric.
\end{lemma}

\begin{proof}
Immediate from Definition~\ref{def:m-compat} and the symmetry of the
compatibility relation on $\Phi_{\geq -1}\,$. 
\end{proof}

We next define $R_m$, the $m$-analogue of~$R$,
and prove an $m$-analogue of Theorem~\ref{th:compat}. 

\begin{definition}
\label{def:Rm}
For $\alpha^k\in\Phi^m_{\geq -1}$, we set 
\[
R_m(\alpha^k)=
\begin{cases}
\ \alpha^{k+1} & \text{if $\alpha\in\Phi_{>0}$ and $k<m$;} \\[.05in]
(R(\alpha))^1 & \text{otherwise.}
\end{cases}
\]
\end{definition}


\begin{theorem}
\label{th:alt def}
The compatibility relation on $\Phi^m_{\geq -1}$ 
has the following properties:
\begin{itemize}
\item[(i) ]
$\alpha^k$ is compatible with $\beta^l$
  if and only if $R_m(\alpha^k)$ is compatible with~$R_m(\beta^l)$; 
\item[(ii) ]
$(-\alpha_i)^1$ is compatible with $\beta^l$ if and only if the simple root
  expansion of $\beta$ does not involve~$\alpha_i$. 
\end{itemize}
Furthermore, conditions {\rm (i)--(ii)} uniquely determine this relation.
\end{theorem}

The proof of Theorem~\ref{th:alt def} is given in Section~\ref{sec:alt
  def}. It does not rely on any statements proved in between. 

\begin{corollary}
\label{cor:tau-switch}
Simultaneously replacing $R=\tau_-\tau_+$ by 
$R=\tau_+\tau_-$ and changing the colors by $i\mapsto m-i+1$ 
does not change the compatibility relation on $\Phi^m_{\geq -1}\,$. 
\end{corollary}

\begin{proof}
Immediate from Theorem~\ref{th:alt def}, since all we are really doing
is replacing $R_m$ by $(R_m)^{-1}$.
\end{proof}

\begin{corollary}
\label{cor:subcomplex}
If $m'\le m$, then $\Delta^{m'}(\Phi)$ is a (vertex-induced) 
subcomplex of $\Delta^{m}(\Phi)$.
\end{corollary}

\begin{proof}
Follows directly from Theorem~\ref{th:alt def}. 
\end{proof}

For a reducible root system~$\Phi=\Phi_1\times\cdots\times\Phi_l$, 
compatibility is defined by analogy with the non-colored case. 
Two colored roots $\alpha^k\in(\Phi_i)^m_{\geq -1}$ and
$\beta^l\in(\Phi_j)^m_{\geq -1}$ are compatible if either $i\neq j$ 
or else $i=j$ and $\alpha^k$ and $\beta^l$ are compatible as elements of
$(\Phi_i)^m_{\geq -1}$. 

For a non-empty subset $J\subset I$, let $\Phi_J$ denote the parabolic
root subsystem of~$\Phi$ spanned by the simple roots $\alpha_j$ for
$j\in J$. 
For $i\in I$, we denote $\br{i}\stackrel{\rm def}{=}I-\set{i}$.
Thus $\Phi_\br{i}$ consists of all roots in $\Phi$ whose simple root
expansion does not involve~$\alpha_i\,$. 

The compatibility relation is preserved under restriction to a parabolic
subsystem: 

\begin{theorem}
\label{th:restrict}
If $\alpha^k,\beta^l\in\Pgem$ and $\alpha,\beta\in\Phi_J\,$, then 
$\alpha$ and $\beta$ are compatible in $\Pgem$ if and only if
they are compatible in $(\Phi_J)^m_{\ge-1}$.
\end{theorem}

The proof of Theorem~\ref{th:restrict} is given in
Section~\ref{sec:restrict} using the classification of finite Coxeter
groups together with the combinatorial models for generalized cluster
complexes of types $ABD$ presented in Section~\ref{sec:classical}.

\begin{definition}
\label{def:gcc}
The \emph{generalized cluster complex}~$\Delta^m(\Phi)$ is the clique complex
for the compatibility relation on $\Phi^m_{\geq -1}$. 
That is, the ground set for $\Delta^m(\Phi)$ is the set of colored
almost positive roots, and a subset of such roots forms a simplex if
and only if any two of them are compatible. 
\end{definition}

If $\Phi$ is irreducible, Corollary~\ref{cor:tau-switch} implies that
the complex $\Delta^m(\Phi)$ is invariant, up to an isomorphism, 
under interchanging $\tau_+$ and~$\tau_-\,$. 
If $\Phi$ is reducible, $\Phi=\Phi_1\times\cdots\times\Phi_l\,$, 
then $\Delta^m(\Phi)$ is, by design, the join of the complexes $\Delta^m(\Phi_j)$.

As in the non-colored case, all maximal simplices of $\Delta^m(\Phi)$
have cardinality~$n$: 

\begin{theorem}
\label{th:pure}
The simplicial complex  $\Delta^m(\Phi)$ is pure of dimension~$n-1$. 
\end{theorem}

\begin{proof}
Let $F$ be a maximal face of $\Delta^m(\Phi)$ and let $\alpha^k$ be a
colored root in $F$.
Choose $i$ so that $R_m^i(\alpha^k)$ is a negative simple root $(\alpha')^1$.
Then $F'\stackrel{\rm def}{=} R_m^i(F)$ is a maximal face of
$\Delta^m(\Phi)$ containing~$(\alpha')^1$. 
By Theorem~\ref{th:alt def}, all other roots in $F'$ are of the form
$\beta^l$ for $\beta\in(\Phi_{\br{i}})_{\ge -1}\,$, and these roots form a maximal
face $F''$ in $\Delta^m(\Phi_{\br{i}})$.
(Here we use Theorem~\ref{th:restrict}.) 
By induction, this simplex is of dimension~$n-2$, so $F$ is of dimension 
$n-1$.
The base of the induction ($n=1$) is trivial;
cf.\ Example~\ref{example:gcc:n=1} below.
\end{proof}

\begin{proposition}
\label{prop:m-fold}
Each codimension~$1$ face of $\Delta^m(\Phi)$ is contained in exactly
$m+1$ faces of maximal dimension. 
\end{proposition}
\begin{proof}
Let $F$ be a codimension~$1$ face of $\Delta^m(\Phi)$ and let $\alpha^k$, 
$i$, $\alpha'$, $F'$ and $F''$ be defined as in the proof of 
Theorem~\ref{th:pure}.
Then $F''$ is a codimension~$1$ face of~$\Delta^m(\Phi_{\br{i}})$, which
by induction is contained in exactly $m+1$ faces of maximal dimension.
Adjoining $(\alpha')^1$ to each of these faces and applying $R_m^{-i}$,
one obtains all maximal faces containing~$F$.
\end{proof}

\begin{prop}
\label{prop:link}
The link of any face in $\Delta^m(\Phi)$ is isomorphic
to a join of generalized cluster complexes of the
form~$\Delta^m(\Phi_J)$,
for some irreducible parabolic root subsystems $\Phi_j\subset\Phi$. 
\end{prop}

\begin{proof}
It is enough to verify the statement for links of vertices.
Applying $R_m$ to a vertex until it becomes 
a negative simple root, the statement follows. 
\end{proof}

\begin{remark}
\label{rem:dual-root-system}
Assume that $\Phi$ is a crystallographic root system, 
and let $\Phi^\vee$ be the root system dual to~$\Phi$.  
As shown in~\cite[Proposition~3.3]{ga}, two roots
$\alpha,\beta\in\Phi_{\geq -1}$ are compatible if and only if the
corresponding coroots $\alpha^\vee$ and $\beta^\vee$ are compatible in
$\Phi^\vee_{\geq -1}$. 
Thus the cluster complex $\Delta(\Phi)$ is canonically isomorphic
to~$\Delta(\Phi^\vee)$. 
It follows directly from Definition~\ref{def:m-compat} that the same
is true about $\Delta^m(\Phi)$ and~$\Delta^m(\Phi^\vee)$. 
In particular, the generalized cluster complexes of types $B_n$
and~$C_n$ are canonically isomorphic. 
In view of this, we do not consider type~$C_n$ in what follows. 
\end{remark}

\begin{remark}
\label{rem:m-cluster-algebras}
It would be very interesting to find a generalization of the notion of
a cluster complex of a cluster algebra of finite type (as defined
in~\cite{ca2}) 
that would yield the complexes~$\Delta^m(\Phi)$. 
This might also lead to an extension of the concept to infinite types. 
\end{remark}

\section{Examples} 
\label{sec:examples}

Here we illustrate the definition of a generalized cluster complex
(see Definition~\ref{def:gcc}) by considering the special cases
where $m\leq 1$ or $n\leq 2$.

\begin{example}[\emph{$m=0$}] 
\label{example:m=0}
The ground set of $\Delta^0(\Phi)$ consists
of the $n$ negative simple roots in~$\Phi$.
Any two such roots are compatible.
Hence $\Delta^0(\Phi)$ is an $(n-1)$-dimensional simplex. 
\end{example}

\begin{example}[\emph{$m=1$}] 
\label{example:m=1}
By design, $\Delta^1(\Phi)=\Delta(\Phi)$ is
the cluster complex of (finite) type~$\Phi$. 
In particular, in type~$A_n$ we obtain Stasheff's associahedron,
while in type $B_n$/$C_n$ we get the Bott-Taubes' cyclohedron. 
See \cite{pcmi,ga} or Section~\ref{sec:classical} below for further
details. 
\end{example}

\begin{example}[\emph{$n=1$}] 
\label{example:gcc:n=1}
For a root system $\Phi$ of type~$A_1$, the simplicial complex
$\Delta^m(\Phi)$ consists of $m+1$ disconnected points
$(-\alpha_1)^1,\alpha_1^1\,,\dots,\alpha_1^m\,$. 
\end{example}

\begin{example}[\emph{$n=2$}] 
\label{example:gcc,n=2} 
Let $\Phi$ be a root system of type~$I_2(a)$. 
The ($1$-dimensional) simplicial complex
$\Delta^m(\Phi)$ is an $(m+1)$-regular graph on $am+2$ vertices. 
Figure~\ref{fig:B2.3} shows this graph for a root system of type
$I_2(4)=B_2$ and $m=3$.
The figure is color-coded to highlight the subgraphs corresponding 
to $m=2$ and $m=1$ and thus to illustrate how $\Delta^m(B_2)$ 
changes with increasing~$m$ (cf.\ Corollary~\ref{cor:subcomplex}).
The simple roots are $\alpha_1$ and~$\alpha_2\,$, with $I_+=\set{1}$
and $I_-=\set{2}$.  The other two
roots are denoted by $\beta_1=2\alpha_1+\alpha_2$ and 
$\beta_2=\alpha_1+\alpha_2$.
The black edges are the edges of $\Delta^1(B_2)$, the black 
and dark-gray edges form $\Delta^2(B_2)$, and the 
edges of all colors form $\Delta^3(B_2)$.

%
%

\begin{figure}[ht]
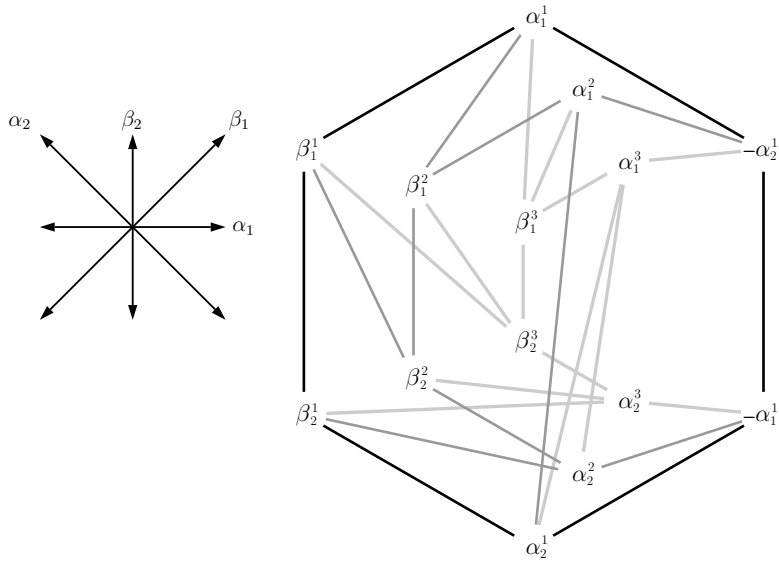

\centerline{
\raisebox{90pt}{
\scalebox{0.7}{
\begin{picture}(13,0)(-63,-50)
                \put(58,-2){$\alpha_1$}
                \put(-63,55){$\alpha_2$}
                \put(56,55){$\beta_1$ 
                            }
                \put(-2,55){$\beta_2$ 
                            }
        \end{picture}
\epsfbox{B2root.ps}
}
}
\hspace{0.2 in}
\scalebox{0.5}{\epsfbox{B2.3.ps}}
}
\caption{Root system of type $B_2$ and the complex $\Delta^3(B_2)$}
\label{fig:B2.3}
\end{figure}

Figure~\ref{fig:B2.3star} shows the same graph, 
drawn so as to make its symmetry apparent.
The map~$R_m$ acts by a $\frac{2\pi}{7}$ counterclockwise rotation of this
picture. 

\begin{figure}[ht]
\centerline{\scalebox{0.5}{\epsfbox{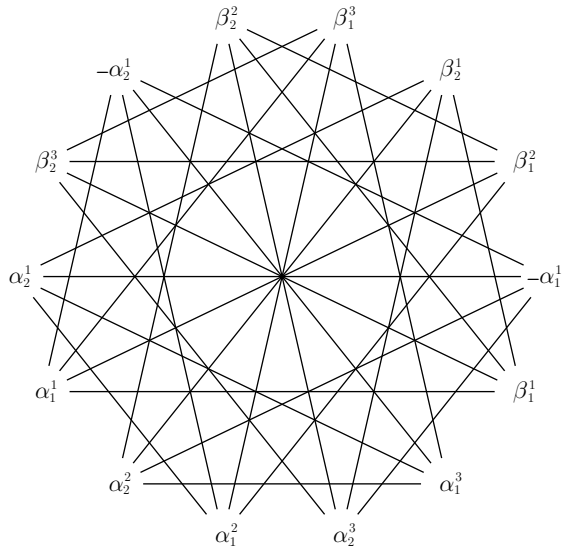}}}
\caption{Representation of $\Delta^3(B_2)$ showing $7$-fold symmetry}
\label{fig:B2.3star}
\end{figure}


In general, the graph $\Delta^m(I_2(a))$ can be constructed in the
plane as follows. 
Take the vertex set to be the integers modulo $(am+2)$, identified
with the vertices of a regular $(am+2)$-gon so that $0,1,2,\dots$ lists them
in clockwise order.
For $a$ odd, the edge set has $(am+2)$-fold rotational symmetry 
and connects each vertex $v$ to the $m+1$ vertices of the form 
$v+\frac{a-1}{2}m+j$, for $j=0,\dots,m$. 
Figure~\ref{fig:I25.4} shows this graph for $a=5$ and $m=4$.
For $a$ even, fix an odd integer~$i$.
The edge set of $\Delta^m(I_2(a))$ has 
$(am/2+1)$-fold rotational symmetry,
and connects $0$ to the vertices $i,i+2,\ldots,i+2m$.
Figure~\ref{fig:B2.3star} is $\Delta^3(B_2)$ drawn in this style with $i=3$.
\end{example}

\begin{figure}[ht]
\centerline{\scalebox{0.5}{\epsfbox{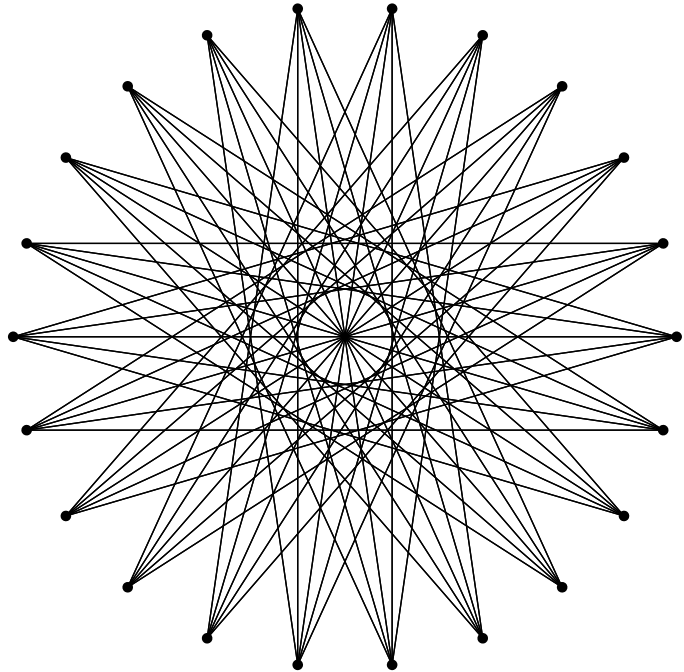}}}
\caption{$\Delta^4(I_2(5))$}
\label{fig:I25.4}
\end{figure}

\begin{figure}[t]
\centerline{\scalebox{0.8}{\epsfbox{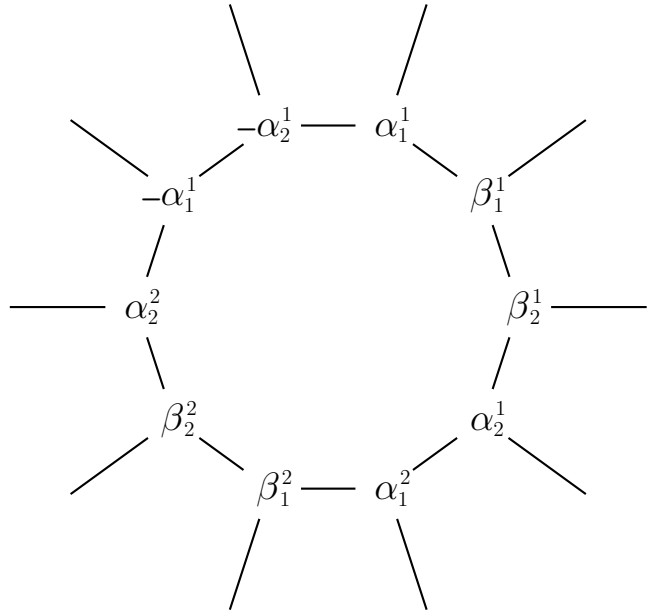}}}
\caption{$\Delta^2(B_2)$ in the real projective plane}
\label{fig:B2.2RP}
\end{figure}

\begin{remark}
\label{rem:RP2}

In the special case $n=m=2$, the complex $\Delta^2(I_2(a))$ 
can be obtained from the boundary of a regular $(2a+2)$-gon
by adding an edge connecting every pair of antipodal vertices,
as illustrated in Figure~\ref{fig:B2.2RP}.
Thus,  $\Delta^2(I_2(a))$ is the \hbox{$1$-skeleton} of a polygonal
subdivision of~$\mathbb{RP}^2$.
It would be interesting to determine which complexes $\Delta^m(\Phi)$
(perhaps all of them?) can be realized as skeleta of polyhedral
$(n+m-2)$-dimensional manifolds. 
\end{remark}

\pagebreak[3]

\section{The classical types}
\label{sec:classical}

In this section, we describe combinatorial models realizing the generalized cluster
complexes $\Delta^m(\Phi)$ for the classical types~$A_n\,$,
$B_n/C_n\,$, and~$D_n\,$.
Each of these models is constructed in terms of collections of
non-intersecting diagonals in a certain convex polygon.
In every case, the isomorphism between the suggested model and the
root-theoretic description of $\Delta^m(\Phi)$ can be verified by
straightforward induction on~$m$.  The base of the induction is the
non-colored ($m=1$) case, described in \cite[Section 3.5]{ga}.

For the types $A_n$ and~$B_n\,$, we rediscover the constructions given earlier
by E.~Tza\-naki~\cite{Tzanaki}. 
The type-$A_n$ construction was inspired by the paper~\cite{przytycki-sikora}
by J.~H.~Przytycki and A.~S.~Sikora. 

\addtocontents{toc}{\SkipTocEntry}
\subsection{Type~$A_n$} 
\label{subsec:gcc-An}

It is easy to see that a convex polygon $\Poly$ with $N$ vertices can be
dissected into (convex) $(m+2)$-gons by pairwise non-crossing
diagonals if and only if $N\equiv 2\bmod m$. 
(Here and in what follows, ``non-crossing'' means that diagonals have
no common points in the \emph{interior} of~$\Poly$.) 
Let $\Poly$ be a convex polygon  with $(n+1)m+2$ vertices. 
To simplify descriptions, we take $\Poly$ to be regular.
A diagonal of~$\Poly$ is called \emph{$m$-allowable} if it
cuts~$\Poly$ into two polygons each of which can be dissected into
$(m+2)$-gons. That is, an $m$-allowable diagonal connects two vertices
such that each of the two paths connecting them along the perimeter
of~$\Poly$ passes through a number of vertices that is divisible by~$m$. 

Let us now define a simplicial complex $\Delta(m,n)$ 
on the set of all $m$-allowable diagonals in~$\Poly$. 
Such diagonals form a simplex in  $\Delta(m,n)$ if and only if they
are pairwise non-crossing.
The complex  $\Delta(m,n)$ is easily seen to be pure, with maximal
simplices formed by $n$-tuples of diagonals that cut $\Poly$ into
$(m+2)$-gons. 

We now construct an isomorphism between $\Delta(m,n)$ 
and the generalized cluster complex~$\Delta^m(\Phi)$  
for a root system $\Phi$ of type~$A_n\,$.
Under the isomorphism, $R_m$ will correspond to 
a clockwise rotation of $\Poly$ taking $P_2$ to $P_1$, etc. 
To avoid additional notation, we refer to this clockwise rotation
as $R_m$ while defining the isomorphism.

\begin{figure}[ht]
\centerline{\scalebox{0.7}{\epsfbox{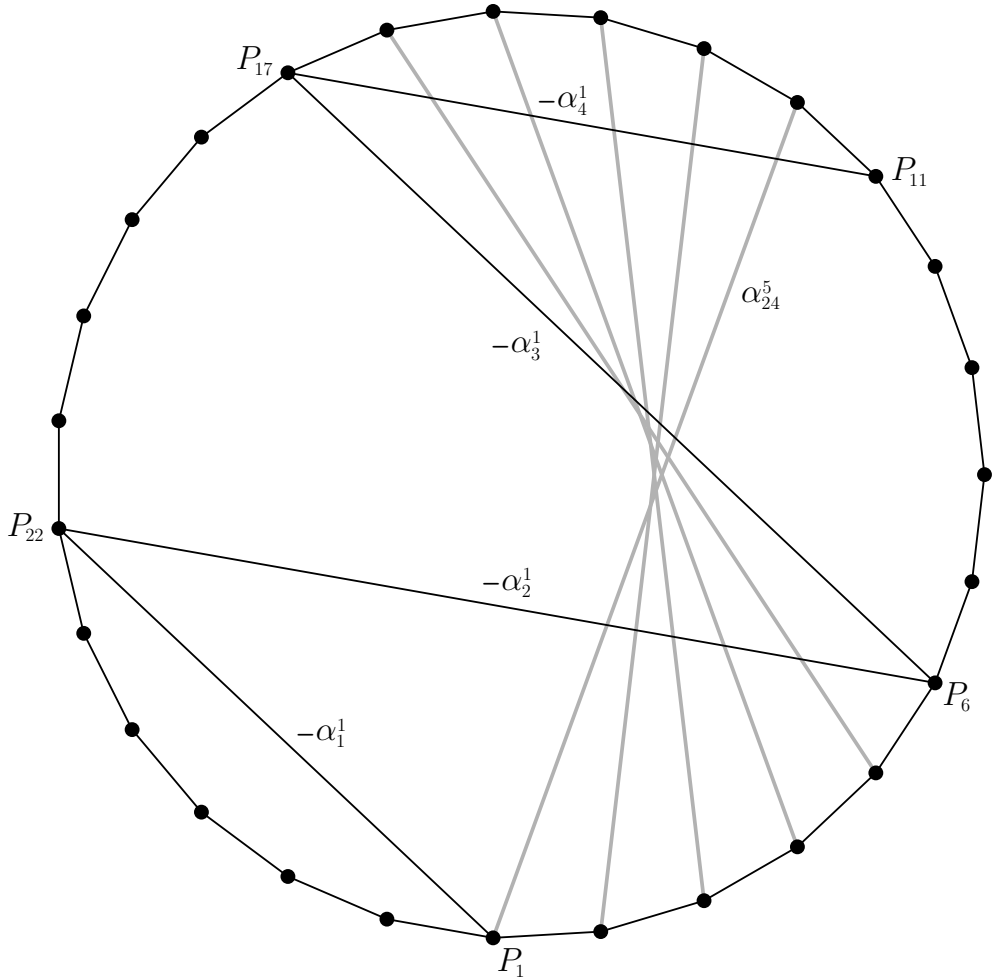}}}
\caption{The $5$-snake for $A_4$}
\label{fig:A4.5snake}
\end{figure}

We use the standard labeling of the simple roots in~$\Phi$
(see, e.g.,~\cite{Bourbaki}) and take $I_+$ to be the odd indices.
Let $P_1,P_2,\ldots,P_{(n+1)m+2}$ be the vertices of $\Poly$, labeled 
counterclockwise.  
For $1\le i\le\frac{n+1}{2}$, identify the negative simple root
$-\alpha_{2i-1}$ with the 
diagonal of $\Poly$ connecting $P_{(i-1)m+1}$ to~$P_{(n+1-i)m+2}$.
For $1\le i\le\frac{n}{2}$, identify $-\alpha_{2i}$ with the 
diagonal connecting $P_{im+1}$ to~$P_{(n+1-i)m+2}$.
Collectively, these $n$ diagonals (each of them $m$-allowable) form
what we call the \emph{$m$-snake} (cf.\ \cite[Figure~3]{ga}). 
The positive roots of $\Phi$ are $\alpha_{ij}=\alpha_i+\cdots+\alpha_j$ for 
each $i$ and $j$ with $1\leq i\leq j\leq n$.
For each such $i$ and $j$,
there are exactly $m$ diagonals which are $m$-allowable and intersect
the diagonals $-\alpha_i,\ldots,-\alpha_j$ and no other diagonals in the $m$-snake.
This collection of diagonals is of the form $R_m^0D,R_m^1D,\ldots,R_m^{m-1}D$
for some diagonal~$D$.
For $1\le k\le m$, we identify $\alpha_{ij}^k$ with $R_m^{k-1}D$.
Figure~\ref{fig:A4.5snake} shows the $m$-snake for $n\!=\!4$ and~$m\!=\!5$, 
along with the diagonals identified with the colored roots 
$\alpha_{24}^1\,,\dots,\alpha_{24}^5\,$. 

For $m=1$, the complex $\Delta(1,n)\cong \Delta(A_n)$ is, by definition
(see, e.g., \cite{pcmi} and references therein), the dual complex for 
(the boundary of) the ordinary $n$-dimensional \emph{associahedron}, 
also known as the \emph{Stasheff polytope}. 
The first non-trivial example with $m>1$ is $\Delta(2,2)\cong \Delta^2(A_2)$, 
the complex of quadrangulations of an octagon. 
(Cf.\ Example~\ref{example:gcc,n=2} with $m=2$ and $a=3$,
and also Remark~\ref{rem:RP2} with $a=3$.) 
This $1$-dimensional complex is a non-planar 3-regular graph
with  $8$ vertices and $12$ edges. 
Thus $\Delta(2,2)$ is homotopy equivalent to a wedge of
5 circles.

\addtocontents{toc}{\SkipTocEntry}
\subsection{Type~$B_n/C_n$}
\label{subsec:gcc-Bn}

When $\Phi$ is of type $B_n$ (or~$C_n$---cf.\
Remark~\ref{rem:dual-root-system}), the complex $\Delta^m(\Phi)$ 
can be realized as follows. 
Let $\Poly$ be a centrally symmetric regular polygon with $2nm+2$ vertices.
The map $R_m$ will correspond to a $\frac{\pi}{nm+1}$ clockwise
rotation about the center of~$\Poly$.
The vertices of the complex are of two kinds:
\begin{itemize}
\item 
the \emph{diameters} of~$\Poly$, i.e., diagonals connecting antipodal
vertices (all such diagonals are $m$-allowable in the sense of
Section~\ref{subsec:gcc-An}); 
\item
the pairs $(D,D')$ of 
distinct $m$-allowable diagonals of $\Poly$ such that $D$ is related to 
$D'$ by a half-turn about the center of $\Poly$.
\end{itemize}
Two vertices are called ``non-crossing'' if no diagonal
representing one vertex crosses a diagonal representing the second vertex.
The simplices in the complex are the sets of pairwise non-crossing vertices.
It is easily verified that this is a pure complex 
of dimension~$n$.
Its maximal simplices correspond to centrally symmetric
dissections of~$\Poly$ into $(m+2)$-gons. 

\begin{figure}[ht]
\centerline{\scalebox{1}{\epsfbox{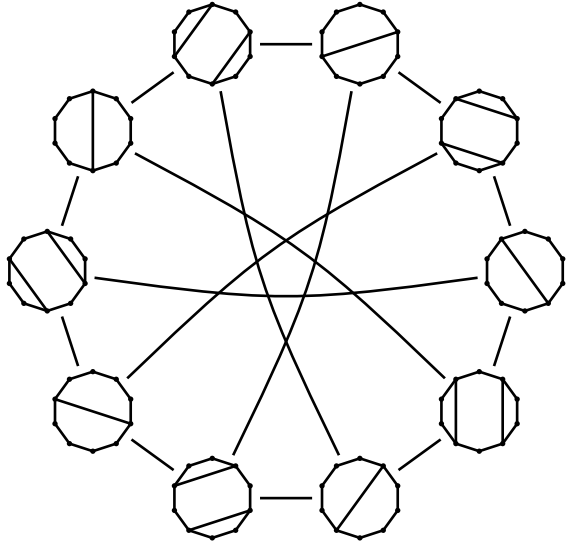}}}
\caption{Combinatorial realization of $\Delta^2(B_2)$}
\label{fig:B2.2}
\end{figure}

Figure~\ref{fig:B2.2} illustrates this complex for $n=m=2$.
The edges of the graph correspond to centrally symmetric
quadrangulations of a $10$-gon.

For $m=1$, this construction specializes to one of the common 
definitions of the $n$-dimensional \emph{cyclohedron},
or \emph{Bott-Taubes polytope} (see, e.g., \cite{pcmi, ga} and
references therein). 

In order to describe an isomorphism between this complex
and~$\Delta^m(B_n)$, we relate the former to $\Delta(m,n)$.
Label the simple roots $\alpha_1,\ldots,\alpha_n$ so that 
$\alpha_n$ is the only simple root in its $W$-orbit. 
The negative simple roots correspond to the orbits, under central symmetry,
of the diagonals on the $m$-snake of type~$A_{2n-1}$ 
(see Section~\ref{subsec:gcc-An}).
If $\beta_1,\ldots,\beta_{2n-1}$ are the simple roots of~$A_{2n-1}\,$, then
$-\alpha_i$ is encoded by the pair of diagonals corresponding to $-\beta_i$ and
$-\beta_{2n-i}\,$, for $1\leq i\leq n-1$, while $-\alpha_n$ is encoded
by the diameter corresponding to~$-\beta_n$.

The positive roots of $B_n$ can be divided into three categories:
\[\begin{array}{rll}
\mbox{(I)}&\alpha_i+\alpha_{i+1}+\cdots+\alpha_j&\mbox{for }i\le j<n\\
\mbox{(II)}&\alpha_i+\alpha_{i+1}+\cdots+\alpha_n&\mbox{for }i\le n\mbox{, and}\\
\mbox{(III)}&\alpha_i+\alpha_{i+1}+\cdots+\alpha_{j-1}+2\alpha_j+2\alpha_{j+1}+\cdots+2\alpha_n&\mbox{for }i<j\le n.
\end{array}\]

A colored positive root $\alpha^k$ for $\alpha$ in category (I) is encoded 
by the pair consisting of the diagonal for $(\beta_i+\beta_{i+1}+\cdots+\beta_j)^k$
and the diagonal for $(\beta_{2n-j}+\beta_{n-j+1}+\cdots+\beta_{2n-i})^k$.
For $\alpha$ in category (II), $\alpha^k$ is encoded by the diameter 
corresponding to $(\beta_i+\beta_{i+1}+\cdots+\beta_{2n-i})^k$.
For $\alpha$ in category (III), $\alpha^k$ is encoded by the diagonal
corresponding to $(\beta_{i}+\beta_{i+1}+\cdots+\beta_{2n-j})^k$ and 
the diagonal corresponding to $(\beta_j+\beta_{j+1}+\cdots+\beta_{2n-i})^k$.

\addtocontents{toc}{\SkipTocEntry}
\subsection{Type~$D_n$}
\label{subsec:gcc-Dn}

Let $\Poly$ be a regular polygon with $2(n-1)m+2$ vertices.
The vertices in the combinatorial realization of~$\Delta^m(D_n)$
fall into two groups. 
The vertices in the first group correspond one-to-one to pairs of distinct 
non-diameter \hbox{$m$-allowable} diagonals in~$\Poly$ related by a half-turn.
In the second group, each vertex is indexed by a diameter of~$\Poly$, 
together with one of two \emph{flavors}, which we will call ``dashed''
and ``gray,'' and picture accordingly. Thus, each diameter 
occurs twice,  in each of the two flavors. 
We label the vertices of $\Poly$ counterclockwise: 
\[
P_1,P_2,\ldots,P_{(n-1)m+1},-P_1,-P_2,\ldots,-P_{(n-1)m+1}. 
\]
We call $[P_1,-P_1]$ the {\em primary diameter}.
By construction, the map $R_m$ acts by rotating $P_2$ clockwise to $P_1$
and switching the flavor of certain diameters.
Specifically, $R_m$ preserves flavor when applied to a diameter
of the form $[P_k,-P_k]$ (for $1\leq k\leq (n-1)m+1$) 
unless $k=1$ or $k\equiv 2 \bmod m$.
If $k=1$ or $k\equiv 2 \bmod m$, the flavor is switched.

By analogy with types $A_n$ and~$B_n\,$, the complex $\Delta^m(D_n)$
is realized as a clique complex for a certain ``compatibility''
relation defined as follows. 
Two vertices at least one of which is not a diameter are compatible or not
according to precisely the same rules as in the type-$B_n$ case.
A more complicated condition determines whether two diameters are
compatible. 
Two diameters with the same location and different flavors are
compatible. 
Two diameters at different locations are compatible if and only if 
applying $R_m$ repeatedly until either of them is in position
$[P_1,-P_1]$ results in diameters of the same flavor.
Explicitly, let $D$ be a flavored diameter and let $\overline{D}$
denote $D$ with its flavor reversed.
Then for $1\leq k\leq (n-1)m$, the flavored diameter $R_m^k(D)$ is
compatible with~$D$ (and incompatible with $\overline{D}$) 
if and only if $\left\lceil\frac{k}{m}\right\rceil$ is even.

\begin{figure}[ht]
\centerline{\scalebox{0.7}{\epsfbox{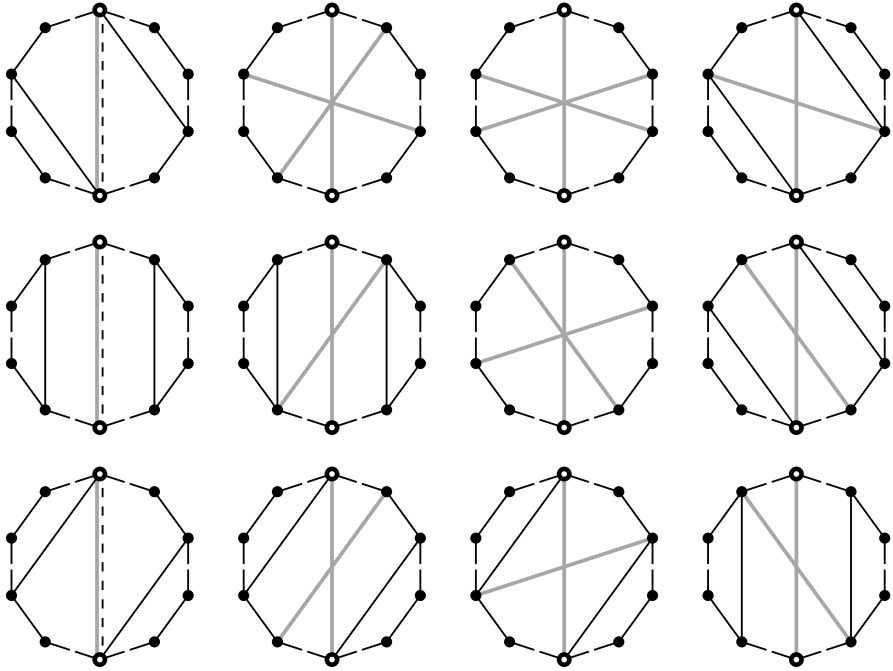}}}
\caption{Some maximal simplices in $\Delta^2(D_3)$}
\label{fig:D3.2diam}
\end{figure}

\begin{figure}[ht]
\centerline{\scalebox{0.7}{\epsfbox{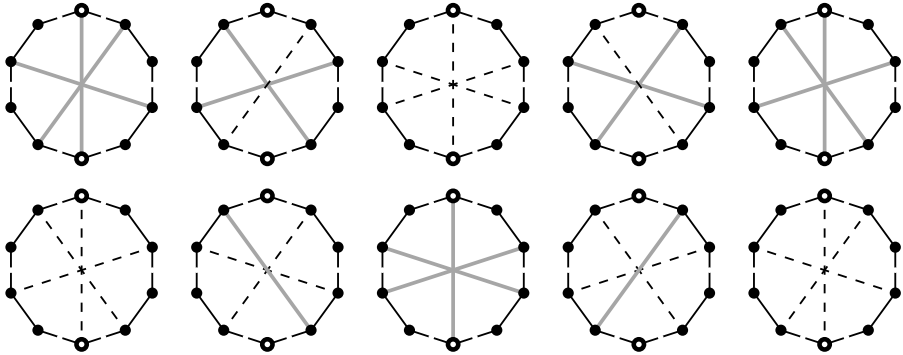}}}
\caption{The orbit of a maximal simplex in $\Delta^2(D_3)$}
\label{fig:D3.2orb}
\end{figure}

Figure~\ref{fig:D3.2diam} shows the $12$ maximal simplices of 
$\Delta^2(D_3)$ containing the gray diameter $[P_1,-P_1]^\text{gr}$.
The points $\pm P_1$ are marked with a white dot.
Some edges of $\Poly$ are broken, indicating the locations at which diameters
change flavors under the action of $R_m$.
Specifically, $P_i$ and $P_{i-1}$ are connected by a broken edge if
and only if $R_m$ maps a flavored diameter $[P_i,-P_i]$ to the opposite-flavored
diameter $[P_{i-1},-P_{i-1}]$.
%
To illustrate the rules of compatibility for pairs of diameters,
Figure~\ref{fig:D3.2orb} shows the orbit
of a maximal simplex under the action of~$R_m$.

We next describe an isomorphism between the complex defined above
and the original construction of~$\Delta^m(D_n)$.
Label the simple roots $\alpha_1,\ldots,\alpha_n$ so that 
$\alpha_{n-1}$ and $\alpha_n$  correspond to 
the ends of the fork in the Dynkin diagram of type~$D_n$.
The negative simple roots correspond to the orbits, under central symmetry,
of the diagonals in the $m$-snake for $A_{2n-3}$, except
that $\alpha_n$ and $\alpha_{n-1}$ both correspond to the same
diameter
(in two different flavors).
Specifically, if $\beta_1,\ldots,\beta_{2n-3}$ are the simple roots of
$A_{2n-3}\,$, then 
$-\alpha_i$ is encoded by the pair of diagonals corresponding to $-\beta_i$ and
$-\beta_{2n-i-2}\,$, for $1\leq i\leq n-2$.
The dashed copy of the diameter corresponding to $-\beta_n$ 
encodes~$-\alpha_{n-1}\,$, while $-\alpha_n$ is encoded by the gray 
diameter corresponding to~$-\beta_n$.
We label the vertices of $\Poly$ so that this is the primary diameter $[P_1,-P_1]$.

The positive roots of $D_n$ can be divided into two categories:
\[\begin{array}{rll}
\mbox{(I)}&\alpha_i+\alpha_{i+1}+\cdots+\alpha_j&\mbox{for }i\le j<n\\
\mbox{(II)}&\alpha_i+\alpha_{i+1}+\cdots+\alpha_{n-2}+\alpha_j+\alpha_{j+1}+\cdots+\alpha_n&\mbox{for }i<j\le n.
\end{array}\]

A colored positive root $\alpha^k$ for $\alpha$ in category (I) is encoded 
by the pair consisting of the diagonal for $(\beta_i+\beta_{i+1}+\cdots+\beta_j)^k$
and the diagonal for $(\beta_{2n-j-2}+\beta_{2n-j-1}+\cdots+\beta_{2n-i-2})^k$
if $j<n-1$.
If $\alpha=\alpha_i+\cdots+\alpha_{n-1}$ then $\alpha^k$ is encoded by 
the gray diameter corresponding to 
$(\beta_i+\beta_{i+1}+\cdots+\beta_{2n-i-2})^k$.

For $\alpha$ in category (II), $\alpha^k$ is encoded by the diagonal
corresponding to $(\beta_{i}+\beta_{i+1}+\cdots+\beta_{2n-j-2})^k$ and 
the diagonal corresponding to $(\beta_j+\beta_{j+1}+\cdots+\beta_{2n-i-2})^k$
if $j<n$.
If $\alpha=\alpha_i+\cdots+\alpha_{n-2}+\alpha_n$ then $\alpha^k$ is encoded by 
the dashed diameter corresponding to 
$(\beta_i+\beta_{i+1}+\cdots+\beta_{2n-i-2})^k$.

\begin{remark}
In types $A_n$ and~$B_n\,$, the combinatorial models for 
the complex $\Delta^m(\Phi)$ presented in Sections~\ref{subsec:gcc-An} 
and~\ref{subsec:gcc-Bn} describe this complex 
as a subcomplex of the ordinary cluster complex $\Delta(\Psi)$ 
for a larger root system~$\Psi$ 
of type $A_{nm+m-1}$ or~$B_{nm}$, respectively. 
It is not clear whether such an embedding exists for the type~$D_n\,$.
\end{remark}

\section{Proof of Theorem~\ref{th:alt def}}
\label{sec:alt def}

The proof makes use of the following results, which can be extracted
from~\cite{ga}. 

\begin{lemma}
{\ }

\begin{enumerate}
\item
Each $R$-orbit in $\Phi_{\geq -1}$ either has size $h+2$ and contains two 
negative simple roots, or has size $(h+2)/2$ and contains
one negative simple root.
\item
If $h$ is even, then for every orbit of size $h+2$, the two
negative simple roots are placed symmetrically in the orbit.
That is, applying $R$ to one of the negative simple roots
$(h+2)/2$ times yields the other negative simple root.
\item
If $h$ is odd, then necessarily all $R$-orbits have size $h+2$.
The two negative simple roots in each orbit are placed so that 
applying $R$ to one of the negative simple roots
$(h+3)/2$ times yields the other negative simple root.
Then applying $R$ an additional $(h+1)/2$ times yields 
the original negative simple root.
\item
In every orbit, positive simple roots are adjacent to their negatives.
When $h$ is odd, both positive simple roots in an orbit 
are located on the longer of the two paths between the two negative
simple roots.
\end{enumerate}
\end{lemma} 

If $m=1$, then both characterizations of compatibility specialize to the
non-colored version, so from now on we assume that $m>1$.
Let $\alpha^k,\beta^l\in\Pgem$. 
Without loss of generality, let $d(\alpha)\le d(\beta)$.
Proving condition (i) requires us to consider several cases.
First, if neither $\alpha$ nor $\beta$ is negative simple 
and neither $k$ nor $l$ equals~$m$, then condition (i) is immediate.

Let us consider the cases where neither $\alpha$ nor $\beta$ is negative 
simple but $k=m$ or $l=m$ or both.
If $k=l=m$, then $\alpha^k$ and $\beta^l$ are compatible if and only if
$\alpha$ and $\beta$ are compatible, and $R_m(\alpha^k)=R(\alpha)^1$
is compatible with $R_m(\beta^l)=R(\beta)^1$ if and only if 
$R(\beta)$ and $R(\alpha)$ are compatible, so (i) follows by the 
$R$-invariance of non-colored compatibility.
If $m=k>l$, then $R_m(\alpha^k)=R(\alpha)^1$ and since 
$d(R(\alpha))<d(\alpha)\le d(\beta)$,
$R(\alpha)^1$ is compatible with $R_m(\beta^l)=\beta^{l+1}$ if and only if 
$R(\alpha)$ is compatible with $\beta$, and this is exactly the requirement
for $\alpha^k$ to be compatible with $\beta^l$.

For $k<l=m$ there are two cases: $d(\alpha)<d(\beta)$ and
$d(\alpha)=d(\beta)$.
If $d(\alpha)<d(\beta)$, then $\alpha^k$ is compatible with 
$\beta^l$ if and only if $\alpha$ is compatible with $\beta$.
Also, $d(\alpha)\le d(R(\beta))$, so
 $R_m(\alpha^k)=\alpha^{k+1}$ is compatible with 
$R_m(\beta^l)=R(\beta)^1$ if and only if $R(\alpha)$ and $R(\beta)$
are compatible, and (i) follows by $R$-invariance of non-colored 
compatibility.
If $d(\alpha)=d(\beta)$, then $\alpha^k$ is compatible with 
$\beta^l$ if and only if $\alpha$ is compatible with $R(\beta)$.
Now $d(\alpha)> d(R(\beta))$, so $R_m(\alpha^k)=\alpha^{k+1}$ 
is compatible with $R_m(\beta^l)=R(\beta)^1$ if and only if 
$\alpha$ and $R(\beta)$ are compatible.

If both $\alpha$ and $\beta$ are negative simple, then $k=l=1$,
$R_m(\alpha^1)=R(\alpha)^1$, and $R_m(\beta^1)=R(\beta)^1$.
Now (i) is immediate by $R$-invariance.

If only one of the two roots is negative simple, then $d(\alpha)\le
d(\beta)$ implies that 
the negative simple root must be~$\alpha$, and in particular $k=1$ and 
$d(\alpha)<d(\beta)$.
The colored roots $\alpha^1$ and $\beta^l$ are compatible if and only if
$\alpha$ and $\beta$ are compatible.
If $l=m$, then $R_m(\alpha^1)=R(\alpha)^1$ is compatible with
$R_m(\beta^m)=R(\beta)^1$ if and only if $R(\alpha)$ and $R(\beta)$
are compatible, and (i) follows by $R$-invariance.
We now assume $l<m$.
If $d(R(\alpha))\ge d(\beta)$, then
$R_m(\alpha^1)=R(\alpha)^1$ is compatible with 
$R_m(\beta^l)=\beta^{l+1}$ if and only if $R(\alpha)$ and 
$R(\beta)$ are compatible, so (i) follows by $R$-invariance.

The condition $d(R(\alpha))\ge d(\beta)$ will hold for most
negative simple roots $\alpha$ and positive roots $\beta$.  
The only way for it to fail is if $h$ is odd, $d(\beta)=(h+1)/2$ 
and $d(R(\alpha))=(h-1)/2$.
In this case, $R_m(\alpha^1)=R(\alpha)^1$ is compatible with 
$R_m(\beta^l)=\beta^{l+1}$ if and only if $R(\alpha)$ is compatible 
with $\beta$.
Thus condition (i) in this case amounts to verifying that $\alpha$ and 
$\beta$ are compatible if and only if $R(\alpha)$ and $\beta$ are compatible.
Because  $d(\beta)=(h+1)/2$, the root $\beta$ lies on the long path between two 
negative simple roots in its $R$-orbit and is adjacent to a negative simple
root.  Therefore $\beta$ must be a (positive) simple root.
If $\beta=-\alpha$, then $\alpha$ and $\beta$ are not compatible, and 
$R(\alpha)=\beta$ so $R(\alpha)$ and $\beta$ are not compatible.
If $\beta\neq-\alpha$, then $\alpha$ and $\beta$ are compatible and 
$R^{-1}(\beta)$ is a negative simple root distinct from $\alpha$, 
so $R^{-1}(\beta)$ and $\alpha$ are compatible.  Therefore
$R(\alpha)$ and $\beta$ are compatible. 

If $\alpha$ is negative simple, then the third condition in 
Definition~\ref{def:m-compat} always applies, and condition
(ii) follows. 

It remains to prove the last statement of Theorem~\ref{th:alt def}.
The following is a generalization of Lemma~\ref{lem:dihedral}. 

\begin{lemma}
\label{lem:m-dihedral}
The order of $R_m$ is $(mh+2)/2$ if 
$w_\circ = -1$, and is $mh+2$ otherwise.
Every $\langle R_m \rangle$-orbit in  $\Phi^m_{\geq - 1}$ has a
nonempty intersection with the set $(- \Pi)^1\stackrel{\rm
  def}{=}\{(-\alpha_i)^1\}_{i\in I}$.
These intersections are precisely the $\langle - w_\circ \rangle$-orbits
in~$(- \Pi)^1$.
\end{lemma}

\begin{proof}
Immediate from Lemma~\ref{lem:dihedral} and Definition~\ref{def:Rm}. 
\end{proof}

It follows from Lemma~\ref{lem:m-dihedral} that conditions (i) and~(ii)
uniquely define the compatibility relation on $\Phi^m_{\geq - 1}\,$.
\qed

\section{Proof of Theorem~\ref{th:restrict}} 
\label{sec:restrict}

\begin{proof}[Type $A_n$]
It is sufficient to treat the case $J=\br{i}$.
Let $\alpha^k\in\Pgem\,$.
By Theorem~\ref{th:alt def}, $\alpha\in\Phi_{\br{i}}$ if and only if 
the $m$-allowable diagonal in $\Poly$ corresponding to $\alpha^k$ 
does not intersect the diagonal $D$ corresponding to $(-\alpha_i)^1$.
The diagonal $D$ dissects $\Poly$ into two polygons $\Poly_{<i}$ and
$\Poly_{>i}$.
The remaining negative simple roots (besides $-\alpha_i$) form two 
$m$-snakes, one in each of the two polygons.

The diagonal corresponding to a colored positive root $\alpha_{ij}^k$  
is uniquely determined as the $k$th (in clockwise order) of the 
exactly $m$ diagonals which intersect the diagonals for 
$-\alpha_i,-\alpha_{i+1},\ldots,\alpha_{j}$ and no other diagonals 
in the $m$-snake.
Thus the correspondence between $(\Phi_{\br{i}})^m_{\ge-1}$ and the 
diagonals in the polygons $\Poly_{<i}$ and $\Poly_{>i}$ is identical 
to the correspondence between colored roots $\alpha^k\in\Pgem$ with 
$\alpha\in\Phi_{\br{i}}$ and diagonals in~$\Poly$.
\end{proof}

\begin{proof}[Type $B_n$]
Theorem~\ref{th:restrict} for $\Phi$ of type $B_n$ follows 
easily from the type-$A_{2n-1}$ case by identifying parabolic  root subsystems 
of $B_n$ with parabolic root subsystems of $A_{2n-1}$ which are fixed by 
the involution $\beta_i\mapsto\beta_{2n-i}$.
\end{proof}

\begin{proof}[Type $D_n$]
Again, it suffices to consider the case $J=\br{i}$.
First, assume that $i\not\in\set{n-1,n}$.
By Theorem~\ref{th:alt def}, $\alpha^k\in\Pgem$ has 
$\alpha\in\Phi_{\br{i}}$ if and only if 
the corresponding diameter or diagonal pair in $\Poly$ does not intersect
the diagonal pair $\{D,D'\}$ corresponding to~$(-\alpha_i)^1$.
These two diagonals dissect $\Poly$ into three polygons, one of which is
centrally symmetric and two of which are related to each other by central 
symmetry.
The proposition follows in this case by an argument analogous to the type $A_n$
case.

By symmetry, we need only consider one additional case: $i=n$.
The colored almost positive roots compatible with $-\alpha_n$ correspond to 
the dashed primary diameter, all gray non-primary diameters, and all
pairs of diagonals which do not intersect the primary diameter.
The parabolic root subsystem $\Phi_{\br{i}}$ is of type $A_{n-1}$.
We will realize $\Delta^m(\Phi_{\br{i}})$ as the complex of non-crossing 
$m$-allowable diagonals in a convex polygon $\Poly'$ with the vertices
\[
P_1,P_2,\ldots,P_{(n-1)m+1},-P_1,Q_1,Q_2,\ldots,Q_m\,,
\]
in counterclockwise order.
The vertices $P_i$ are vertices of $\Poly$ and the vertices $Q_i$ are new.

For a diagonal pair $\{D,D'\}$ in $\Poly$ which does not intersect the
primary diameter, exactly one of the diagonals $D$ and $D'$ connects vertices in 
$\{P_1,P_2,\ldots,P_{(n-1)m+1},-P_1\}$, and is therefore a diagonal in~$\Poly'$.
The dashed primary diameter of $\Poly$ corresponds to the diagonal $[P_1,-P_1]$
of $\Poly'$.
The correspondence between gray diameters in $\Poly$ and the remaining diagonals of
$\Poly'$ is as follows.
The gray diameter $[P_i,-P_i]$ in $\Poly$ is mapped to the diagonal $[P_i,Q_{k(i)}]$
in $\Poly'$ for $k(i)\equiv (i-1)\bmod m$.
Figure~\ref{fig:diag map} illustrates this correspondence 
(for $m=3$ and $\Phi$ of type~$D_5$) by showing a maximal 
simplex in $\Delta^m(\Phi_{\br{i}})$ represented as a collection of gray diameters
in $\Poly$ and as the corresponding collection of diagonals in $\Poly'$.

\begin{figure}[ht]
\[
\begin{array}{ccc}
\scalebox{0.8}{\epsfbox{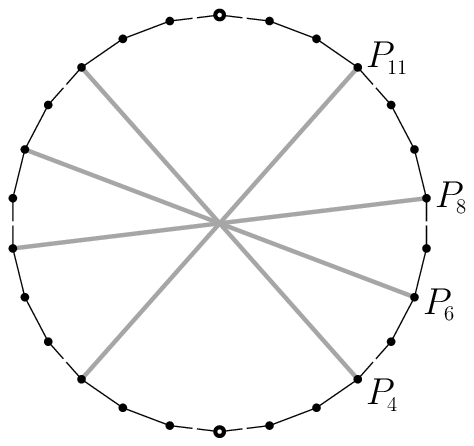}}& 
&\scalebox{0.8}{\epsfbox{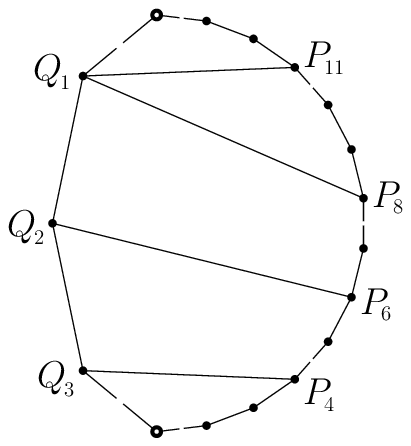}}\\
\mbox{\large$\Poly$}&&\mbox{\large$\Poly'$}
\end{array}
\]
\caption{}
\label{fig:diag map}
\end{figure}

It is immediate that if two diagonal pairs in $\Poly$ do not intersect the primary
diameter, then they intersect each other if and only if the corresponding 
diagonals in $\Poly'$ intersect each other.
It is also immediate that a diagonal pair in $\Poly$ that does not
intersect the primary diameter intersects a gray diameter in $\Poly$ if and only if 
the corresponding diagonals in $\Poly'$ intersect each other.
Intersections with the dashed primary diameter of $\Poly$ also correspond to 
intersections with the diagonal $[P_1,-P_1]$ in $\Poly'$.
One can check that for $1<i<j<(n-1)m+1$,
a gray diameter $[P_i,-P_i]$ is compatible with a gray diameter $[P_j,-P_j]$ 
if and only if the fractional part of $\frac{i-2}{m}$ is greater than or equal to
the fractional part of $\frac{j-2}{m}$.
This is exactly the requirement that $[P_i,Q_{k(i)}]$ and
$[P_i,Q_{k(j)}]$ do not intersect.
Thus our correspondence also maps the compatibility relation among
diameters and $m$-allowable diagonal pairs in $\Poly$ compatible with 
the gray primary diameter to the compatibility relation among $m$-allowable 
diagonals of $\Poly'$.

It remains to verify that the process of mapping from
$(\Phi_{\br{i}})^m_{\ge-1}$ to diagonal pairs and flavored diameters
in~$\Poly$, then mapping to diagonals in $\Poly'$ produces 
the same result as the process of mapping directly from
$(\Phi_{\br{i}})^m_{\ge-1}$ to diagonals 
in $\Poly'$ by the type $A_{n-1}$ construction.
This is a straightforward check. 
\end{proof}

\begin{proof}[Exceptional types]
For all exceptional types (including the non-crystallo\-graphic ones),
Theorem~\ref{th:restrict} was checked with the help of a
computer.
In Definition~\ref{def:m-compat}, the compatibility of $\alpha^k$ and $\beta^l$ depends on the comparison between $k$ and $l$ but not on the actual values of $k$ and $l$.
Thus it is sufficient to check the theorem in the case $m=2$, so that in particular, for each type the theorem reduces to a finite computation.
\end{proof} 

\begin{proof}[Dihedral types $I_2(a)$]
In this case, the statement is trivial
(cf.\ Examples~\ref{example:gcc:n=1}--\ref{example:gcc,n=2}).
\end{proof} 

\section{Face numbers}
\label{sec:face-numbers}

In this section, we enumerate the faces of various dimensions in the
generalized cluster complexes~$\Delta^m(\Phi)$.
Let $f_k(\Phi,m)$ denote the $k$th \emph{face number}~$\Delta^m(\Phi)$ or,
more precisely, the number of $k$-element 
simplices
(i.e., $(k\!-\!1)$-dimensional faces) in~$\Delta^m(\Phi)$. 
We denote by
$N(\Phi,m)=f_n(\Phi,m)$ 
the number of maximal simplices in~$\Delta^m(\Phi)$. \linebreak[3]
We note that $f_0(\Phi)=1$ for any~$\Phi$. 

\begin{example}[\emph{$n=1$}] 
\label{example:f:n=1}
In view of Example~\ref{example:gcc:n=1}, we have
\begin{align*}
f_1(A_1,m) &= m+1. 
\end{align*}
\end{example}

\begin{example}[\emph{$n=2$}] 
\label{example:f:n=2}
By Example~\ref{example:gcc,n=2}, 
the complex $\Delta^m(I_2(a))$ is an $(m+1)$-regular graph on $am+2$ vertices.
Hence 
\begin{align*}
f_1(I_2(a),m) &= am+2, \\
f_2(I_2(a),m) &= \frac{(am+2)(m+1)}{2}. 
\end{align*}
\end{example}



The face numbers $f_k(\Phi,m)$ satisfy the following recurrences which generalize
\cite[Proposition~3.7]{ga}.

\begin{proposition}
\label{prop:fk-recur}
Let $k$ be a positive integer. 
For a reducible root system $\Phi=\Phi_1\times \Phi_2$,
\begin{equation}
\label{eq:fk-reducible}
f_k(\Phi_1\times \Phi_2,m)=\sum_{k_1+k_2=k} f_{k_1}(\Phi_1,m)\,
f_{k_2}(\Phi_2,m). 
\end{equation}
In particular,
\[
N(\Phi_1\times \Phi_2,m)=N(\Phi_1,m)\,N(\Phi_2,m).
\]

\pagebreak[3]

If $\Phi$ is irreducible, then
\begin{equation}
\label{eq:recur-fk}
f_k(\Phi,m) = 
\frac{mh+2}{2k}\sum_{i\in I} f_{k-1}(\Phi_\br{i},m). 
\end{equation}
In particular,
\begin{equation}
\label{eq:recurr-N} 
N(\Phi,m)=\frac{mh+2}{2n}\sum_{i\in I}N(\Phi_{\br{i}},m).
\end{equation} 
\end{proposition}

\begin{proof}
Formula \eqref{eq:fk-reducible} follows from the fact that
$\Delta^m(\Phi_1\times\Phi_2)$ is a join of $\Delta^m(\Phi_1)$ and~$\Delta^m(\Phi_2)$. 
The proof of \eqref{eq:recur-fk} is analogous to the proof of
\cite[Proposition~3.7]{ga}. 
Let us count in two different ways the number of pairs $(\alpha^j,S)$
where $S$ is a $k$-element simplex in $\Delta^m(\Phi)$ and
$\alpha^j\in S$.
On the one hand, the number of such pairs is $k\cdot f_k(\Phi,m)$.
On the other hand, combining Lemma~\ref{lem:m-dihedral}, Theorem~\ref{th:alt def},
and Theorem~\ref{th:restrict}, we conclude that the colored roots $\alpha^j$
belonging to each $\langle R_m\rangle$-orbit
$\Omega\subset\Phi^m_{\geq -1}$ contribute
\[
\frac{mh+2}{2}\sum_{-\alpha_i\in\Omega} f_{k_1}(\Phi_{\br{i}}\,,m) 
\]
to the count (per negative simple root $(-\alpha)^1$), 
implying~\eqref{eq:recur-fk}. 
\end{proof}


It follows from \eqref{eq:recur-fk} by induction on~$m$ that
$f_k(\Phi,m)$ is a polynomial in~$m$ of degree~$k$. 
In this section, we calculate these polynomials for all irreducible
root systems~$\Phi$. 
(The formulas for reducible cases follow at once.) 
It turns out that, with a single exception of $f_4(E_8,m)$,  
the face numbers $f_k(\Phi,m)$ factor over the rationals into linear
factors. 

We present the special case $k=n$ first. 

For the rest of this section, we assume that the root system~$\Phi$ is
irreducible. 
Let $e_1\,,\dots,e_n$ be the \emph{exponents} of~$\Phi$. 

\begin{proposition}
\label{pr:product-N}
The number of maximal simplices in
$\Delta^m(\Phi)$ is given by 
\begin{equation}
\label{eq:product-N}
N(\Phi,m)=\prod_{i}\frac{mh+e_i+1}{e_i+1}.
\end{equation}
\end{proposition}

Proposition~\ref{pr:product-N} is a special case of the general
formula for the face numbers $f_k(\Phi,m)$ given in Theorem~\ref{th:product-f}
below. 

The product appearing in the right-hand side of \eqref{eq:product-N} 
comes up in a variety of contexts related to the combinatorics of the
root system~$\Phi$, the corresponding Coxeter group~$W$, and associated
Lie algebras and Lie groups. 
We refer the reader to \cite{Ath} and references therein;  
see also the discussion in Section~\ref{sec:h-vectors} below. 


One would hope for a uniform formula for the face numbers
in the spirit of~\eqref{eq:product-N}.  This does not quite happen.  
The factorization of $f_k(\Phi,m)$ 
into irreducible polynomials in~$m$ yields a subset of the factors appearing
on the right-hand side of~\eqref{eq:product-N} together with,
in types $D_n$ ($n\geq 4$), $E_6\,$, $E_7\,$, $E_8\,$, $F_4\,$, and~$H_4\,$,
one additional ``mysterious'' factor, which we denote by $c_f(\Phi,k,m)$. 
We emphasize that $c_f(\Phi,k,m)$ is
not a ``fudge factor'' needed to force our formulas to look
like~\eqref{eq:product-N}.  Rather, it is forced upon us by the unique
factorizations of the polynomials $f_k(\Phi,m)$ into irreducibles.

\pagebreak[3]

\begin{theorem}
\label{th:product-f}
The face numbers of the generalized cluster complex $\Delta^m(\Phi)$ 
for an irreducible root system~$\Phi$ 
are given by 
\begin{equation}
\label{eq:fk-product}
f_k(\Phi,m) 
= c_f(\Phi,k,m) \binom{n}{k} \prod_{\operatorname{level}(e_i)\leq k} \frac{mh+e_i+1}{e_i+1}  
 \,. 
\end{equation} 
Here, the ``levels'' of exponents $e_i$ are defined in Table~\ref{table:levels},
with each column showing the exponents of appropriate type,
each of them placed in the row corresponding to its level. 
The factor $c_f(\Phi,k,m)$ is a polynomial in $m$ given as follows:
\begin{itemize}
\item
in types $A_n$, $B_n$, $H_3$, and $I_2(m)$, 
we have $c_f(\Phi,k,m)=1$;
\item
in type~$D_n$, 
we have:  
\[
c_f(D_n,k,m)=
\begin{cases}
1 & \text{if $k\in\{0,1,n-1,n\}$;}\\
\displaystyle\frac{n(n-1)m+k(k-1)m+kn}{kn} & \text{otherwise};
\end{cases}
\]
\item
in types $E_6$, $E_7$, $E_8$, $F_4$, and $H_4$, the factors
$c_f(\Phi,k,m)$ 
are shown in Table~\ref{table:correction_f}. 
\end{itemize}
In particular, for the classical types $A_n$, $B_n$, and $D_n$, the
face numbers are given by 
\begin{align}
\label{eq:fk-An}
f_k(A_n,m)&=\frac{1}{k+1}\binom{n}{k}\binom{(n+1)m+k+1}{k}, \\
\label{eq:fk-Bn}
f_k(B_n,m)&=\binom{n}{k}\binom{nm+k}{k}, \\
\label{eq:fk-Dn}
f_k(D_n,m)&=\binom{n}{k}\binom{(n-1)m+k}{k}+\binom{n-2}{k-2}\binom{(n-1)m+k-1}{k}.
\end{align}
\end{theorem}

\begin{table}[ht]
\begin{center}
\begin{tabular}{|c||c|c|c|c|c|c|c|c|c|c|}
\hline
level & $A_n$  & $B_n$    & $D_n$    & $E_6$ & $E_7$   & $E_8$   & $F_4$ & $H_3$& $H_4$ & $I_2(a)$ \\ \hline
&&&&&&&&&&\\[-3.5mm]
\hline
$1$ & $1$      & $1$      & $1$      & $1$   & $1$     & $1$     & $1$   & $1$ & $1$ & $1$ \\\hline
$2$ & $2$      & $3$      &          &       &         &         &       & $5$ &  & $a-1$ \\\hline
$3$ & $3$      & $5$      & $3$      & $4$   & $5$     & $7$     & $5,7$ & $9$ & $11,19$ &  \\\hline
$4$ & $4$      & $7$      & $5$      & $5$   & $7$     &         & $11$  &  & $29$ &  \\\hline
$5$ & $5$      & $9$      & $7$      & $7,8$ & $9$     & $11,13$ &       &  &  &  \\\hline
$6$ &  $6$ &  $11$&  $9$& $11$  & $11,13$ & $17$    &       &  &  &  \\\hline
$7$ &  $7$ &  $13$&  $11$&       & $17$    & $19,23$ &       &  &  &  \\\hline
$8$ &  $8$ &  $15$&  $13$&       &         & $29$    &       &  &  &  \\\hline
$\vdots$ &  $\vdots$ &$\vdots$   & $\vdots$ & &&&&&& \\\hline
$n\!-\!1$ & $n\!-\!1$  & $2n\!-\!3$ & $2n\!-\!5,\,n\!-\!1$ & & & & & & & \\\hline
$n$ & $n$      & $2n\!-\!1$ & $2n\!-\!3$ & & & & & & &  \\\hline
\end{tabular}

\bigskip

\end{center}
\caption{Levels of exponents}
\label{table:levels}
\end{table}

\vspace{-.2in}

\newlength{\tableskip}
\newlength{\negtableskip}
\setlength{\tableskip}{1.0mm}
\setlength{\negtableskip}{-4.0mm}

\begin{table}[ht]
\begin{center}
\begin{tabular}{|c||c|c|c|c|c|c|}
\hline
$k$ & $D_8$ & $E_6$ & $E_7$   & $E_8$   & $F_4$ & $H_4$\\ \hline
&&&&&&\\[-3.5mm]
\hline
0 & $1$&$1$&$1$&$1$&$1$&$1$\\
\hline
1 & $1$&$1$&$1$&$1$&$1$&$1$\\ 
\hline
&&&&&&\\[\negtableskip]
2 & $\frac{29m}{8}+1$&$\frac{14m}{5}+1$&$\frac{7m}{2}+1$&$\frac{35m}{8}+1$&$\frac{13m}{6}+1$&$\frac{31m}{12}+1$\\[\tableskip]
\hline
&&&&&&\\[\negtableskip]
3 & $\frac{31m}{12}+1$&$\frac{9m}{4}+1$&$\frac{27m}{10}+1$&$\frac{45m}{14}+1$&$1$&$1$\\[\tableskip]
\hline
&&&&&&\\[\negtableskip]
4 & $\frac{17m}{8}+1$&$\frac{5m}{3}+1$&$\frac{21m}{10}+1$&$\frac{46m^2}{7}+\frac{179m}{35}+1$&$1$&$1$\\[\tableskip]
\hline
&&&&&&\\[\negtableskip]
5 & $\frac{19m}{10}+1$&$1$&$\frac{23m}{14}+1$&$2m+1$&&\\[\tableskip]
\hline
&&&&&&\\[\negtableskip]
6 & $\frac{43m}{24}+1$&$1$&$1$&$\frac{13m}{8}+1$&&\\[\tableskip]
\hline
7 & $1$&&$1$&$1$&&\\
\hline
8 & $1$&&&$1$&&\\
\hline
\end{tabular}
\end{center}

\bigskip

\caption{The factors $c_f(\Phi,k,m)$}
\label{table:correction_f}
\end{table}

The type~$A_n$ case of Theorem~\ref{th:product-f} 
is due to J.~H.~Przytycki and
A.~S.~Sikora~\cite{przytycki-sikora}. 
The type $B_n$ case is due to E.~Tzanaki~\cite{Tzanaki}.
The formula \eqref{eq:fk-Dn} appeared in \cite{Tzanaki-1}
in connection with a simplicial complex non-isomorphic to~$\Delta^m(D_n)$.

For comparison, Table~\ref{table:correction_f}
also includes the factors $c_f(\Phi,k,m)$ for the type~$D_8$. 

\pagebreak[3]

The constant term of the factor $c_f(\Phi,k,m)$ as a
polynomial in~$m$ is always~$1$. 
The degree of $c_f(\Phi,k,m)$ is equal to $k$ 
minus the number of exponents at levels~$\leq k$. 
This degree is $\leq 1$ (i.e., $c_f(\Phi,k,m)$ equals~$1$ or
is linear) in every case except for $c_f(E_8,4,m)$, which is quadratic in~$m$. 


We call the face numbers $f_k(A_n,m)$ 
for the generalized cluster complex of type~$A_n$ 
the \emph{Przytycki-Sikora numbers}. 
Thus, $f_k(A_n,m)$ is the number of $k$-tuples of pairwise-noncrossing
diagonals that dissect a given convex $((n+1)m+2)$-gon~$\Poly$ 
into smaller $((l+1)m+2)$-gons ($0\leq l\leq n$). 

The Przytycki-Sikora numbers \eqref{eq:fk-An} 
generalize both the \emph{Fuss numbers} 
\begin{equation}
\label{eq:cluster-fuss}
N(A_n,m)=\frac{1}{n+1}\binom{(n+1)(m+1)}{n}
=\frac{1}{nm+m+1}\binom{(n+1)(m+1)}{n+1}\,, 
\end{equation}
which count the maximal simplices in $\Delta^m(A_n)$,
and the \emph{Kirkman-Cayley numbers} 
\begin{equation}
\label{eq:cluster-cayley}
f_k(A_n,1)=
\frac{1}{k+1}\binom{n}{k}\binom{n+k+2}{k}\,, 
\end{equation}
the face numbers of an ordinary associahedron. 

Both the Fuss numbers and the Kirkman-Cayley numbers 
(each of which in turn generalize the Catalan numbers) 
come up in many algebraic and combinatorial contexts;
see, for example, \cite{przytycki-sikora, sloane-A001764, stanley-dissections}.
For example, the Kirkman-Cayley numbers count plane trees with 
$n+k+3$ vertices, $n+2$ of which are leaves, such that no vertex has 
exactly one successor \cite[Exercise 6.33.c]{EC2}.
The Fuss numbers: 
\begin{itemize}
\item
are dimensions of the Fuss-Catalan algebras 
of V.~Jones and D.~Bisch~\cite{bisch-jones}; 
\item
count noncrossing partitions into blocks of size~$m+1$ 
(see, e.g.,~\cite{przytycki-sikora}); 
\item
count paths on the coordinate plane which connect the points $(0,0)$ and
$((m+1)n,0)$, have steps of the form $(1,1)$ and~$(1,-m)$,  
and never dive below the $x$-axis (see, e.g.,~\cite{GKP}); 
\item
count plane $(m+1)$-ary trees with $(n+1)(m+1)+1$ vertices
(see, e.g.,~\cite[Proposition 6.2.1]{EC2});
\item
count plane rooted trees with $(m-1)[(n+1)(m+1)+1]$ black vertices 
and $(n+1)(m+1)+1$ white vertices such that every black vertex 
is a leaf and every white vertex has $m-1$ black successors (see, e.g.,~\cite{pak}).
\end{itemize} 

The Przytycki-Sikora numbers have interpretations in most of these
contexts. 
It would be interesting to find analogues of these interpretations
for other classical types (or even arbitrary finite types). 

We conclude this section by a curious observation. 
The following result appeared in \cite[Corollaries 3.2
  and~3.4]{Tzanaki}, with a different proof. 

\begin{corollary} 
\label{cor:number-diam}
In the combinatorial model of the complex 
$\Delta^m(B_n)$
presented in Section~\ref{subsec:gcc-Bn}, 
the number of $k$-element simplices that include a diameter is
equal to~$\binom{n-1}{k-1} \binom{nm+k}{k}$.
\end{corollary}

\begin{proof}
The number of diameters in $\Poly$ is $nm+1$. 
By formula~\eqref{eq:fk-An}, the number of $k$-simplices in
$\Delta^m(B_n)$ which include a diameter is equal to 
\[
(nm+1) f_{k-1}(A_{n-1},m) 
=\frac{nm+1}{k}\binom{n-1}{k-1}\binom{nm+k}{k-1}
=\binom{n-1}{k-1} \binom{nm+k}{k}, 
\]
as claimed. 
\end{proof}

Since 
\[
\binom{n-1}{k-1} \binom{nm+k}{k}
=\frac{k}{n}\binom{n}{k} \binom{nm+k}{k}
=\frac{k}{n} f_k(B_n,m),
\]
Corollary~\ref{cor:number-diam} can be restated as saying that the 
proportion of $k$-element simplices in $\Delta^m(B_n)$ that
include a diameter is equal to~$\frac{k}{n}$. 
The latter statement has the following probabilistic
reformulation. 

\begin{corollary} 
\label{cor:prob-diam}
Fix~$k$.
Choose a $k$-element simplex in $\Delta^m(B_n)$ uniformly at random. 
It consists of $k$ pairs of centrally symmetric (possibly identical)
diagonals. 
Pick one of these $k$ uniformly at random. 
Then the probability that it is a diameter is equal~to~$\frac{1}{n}$.
\end{corollary}

It would be interesting to find a direct probabilistic proof of
Corollary~\ref{cor:prob-diam}, leading to a simpler proof of~\eqref{eq:fk-Bn}.

\section{Proof of Theorem~\ref{th:product-f}}
\label{sec:proof-product-f}

\begin{proof}[Exceptional types $E_6$, $E_7$, $F_4$, $G_2$, $H_3$, and~$H_4$]
For each of these types, Theorem~\ref{th:product-f} 
has been verified on a computer by calculating the face numbers $f_k(\Phi,m)$
using the recurrences in Proposition~\ref{prop:fk-recur}. 
\end{proof}

\begin{proof}[Type~$I_2(a)$]
Immediate from Example~\ref{example:f:n=2}. 
\end{proof}

\begin{proof}[Type~$A_n$]
In light of the isomorphism between $\Delta(m,n)$ and $\Delta^m(A_n)$
(see Section~\ref{subsec:gcc-An}),
formula \eqref{eq:fk-An} is equivalent to \cite[Corollary~2]{przytycki-sikora}.
\end{proof}

\begin{proof}[Type~$B_n$]
In view of the combinatorial model for $\Delta^m(B_n)$ presented in
Section~\ref{subsec:gcc-Bn}, 
formula \eqref{eq:fk-Bn} is equivalent to \cite[Theorem
  1.2]{Tzanaki}. 
\end{proof}

The proofs for the types $A_n$ and $B_n$ given in
\cite{przytycki-sikora} and \cite{Tzanaki} are combinatorial 
(bijective). 
Rather than deducing \eqref{eq:fk-An} and \eqref{eq:fk-Bn} directly
from the recurrences in Proposition~\ref{prop:fk-recur}, 
those proofs proceed by establishing bijections with combinatorial
objects that are easily enumerated by the right-hand-sides of
\eqref{eq:fk-An} and \eqref{eq:fk-Bn}. 
Our proof for the type~$D_n$ presented below is combinatorial as well. 

\begin{proof}[Type~$D_n$]
This is the hardest case. 
The proof given here relies in large part on ideas developed by
E.~Tzanaki in an earlier analysis \cite{Tzanaki-1} of an alternative
type-$D$ analogue of~$\Delta(m,n)$. 

\pagebreak[3]

We will work with the combinatorial model of $\Delta^m(D_n)$ described in
Section~\ref{subsec:gcc-Dn}. 
We first define, for each diameter or diagonal pair, a {\em
starting vertex}.  For a diameter, the starting vertex is always the positive vertex 
(i.e., $P_l$ as opposed to $-P_l$).  A diagonal pair has four vertices, and 
the starting vertex is the unique vertex which is positive 
and from which one can travel along a diagonal in the pair while
keeping the center of $\Poly$ on the left.
We will need the following lemma about the maximal faces of 
$\Delta^m(D_n)$ which consist entirely of diameters.

\begin{lemma}
\label{lem:all diams}
For $n>2$, given $n$ diameters of $\Poly$, the following are equivalent.
\begin{enumerate}
\item[(i)]
There exists a choice of flavor for each diameter such that
the flavored diameters are pairwise compatible.
\item[(ii)]
If $1\le a_1\le a_2\le\cdots\le a_n\le (n-1)m+1$ are the indices of the
starting vertices of the diameters, then for each $1\leq j\leq n$, 
we have $a_{j+1}-a_j\le m$ 
(with the convention $a_{n+1}=a_1+(n-1)m+1$).
\end{enumerate}
If these conditions hold, then there are exactly two ways to assign
flavors to the diameters so that they are pairwise compatible.
These two ways of assigning flavors are related by
switching the flavor of each of the $n$ diameters.
\end{lemma}

\pagebreak[3]

\begin{proof}
First consider the case $a_1=1$, so that
one of these diameters is primary.
All of the diameters must be of the same flavor if they are all to be 
compatible with a diameter in primary position.
Take any two of the diameters not in primary position and let $i$ and $k$
be the indices of their starting vertices, with $i<k$.
Then the two (same-flavored) diameters are compatible 
 if and only if, when $[P_i,-P_i]$
is moved to the primary position by~$R_m\,$, and $R_m$ is applied to 
$[P_j,-P_j]$ the same number of times, then both of them experience the same
number of color changes ($\bmod\, 2$). \linebreak[3]
As noted in Section~\ref{sec:restrict}, this occurs if and only if the
fractional part of $\frac{i-2}{m}$ 
is greater than or equal to the fractional part of $\frac{j-2}{m}$.
Thus the $n$ diameters are pairwise compatible (when flavored identically)
if and only if
\[
   \mbox{frac.\ part}\left(\frac{a_2-2}{m}\right)
\ge\mbox{frac.\ part}\left(\frac{a_3-2}{m}\right)
\ge\cdots
\ge\mbox{frac.\ part}\left(\frac{a_n-2}{m}\right).
\]
Since $a_1=1$, this is now easily checked to be equivalent to~(ii).
Furthermore, there are two
choices of flavorings, either all gray or all dashed.

If $a_1>1$, then rotate the set of diameters clockwise until $a_1=1$ 
and apply the argument above.
If the rotated diagonals are pairwise compatible when flavored identically,
apply $R_m^{-1}$ repeatedly to undo the rotation.
This will define a flavoring of the $n$ diagonals such that they are
pairwise compatible.
\end{proof}

We are now prepared to complete the proof of~\eqref{eq:fk-Dn}. 
Given a face $F$ in $\Delta^m(D_n)$ containing at least one diameter, call $F$
a {\em gray face} if for the minimal $i$ such that $F$ contains a diameter in 
position $[P_i,-P_i]$, there is a gray diameter in that position.
Similarly define a {\em dashed face}. A face with two
coinciding diameters is both gray and dashed.

Following~\cite{Tzanaki-1}, we rewrite \eqref{eq:fk-Dn} as the sum of the expressions
\begin{eqnarray}
\label{piece1}
\binom{n-1}{k}\binom{(n-1)m+k}{k}&+&\binom{n-2}{k-2}\binom{(n-1)m+k}{k},\mbox{ and}\\[.03 in]
\label{piece2}
2\binom{n-2}{k-2}\binom{(n-1)m+k}{k}&-&\binom{n-2}{k-2}\binom{(n-1)m+k-1}{k-1},
\end{eqnarray}
and identify the terms of these expressions.
The first term of (\ref{piece1}) is $f_k(B_{n-1},m)$, which also counts the 
$k$-element faces of
$\Delta^m(D_n)$ with at most one gray diameter and no dashed diameters.
By~\cite[Corollary 3.4]{Tzanaki}, the second term is the number of faces of 
$\Delta^m(B_{n-1})$ with $k$ elements, exactly one of which is a diameter.
This also counts the $k$-element faces of $\Delta^m(D_n)$ with exactly
one dashed diameter 
and no gray diameters.
Thus (\ref{piece1}) counts faces of $\Delta^m(D_n)$ with at most one diameter.
Therefore (\ref{piece2}) counts faces of $\Delta^m(D_n)$ with two or more diameters.
The second term of (\ref{piece2}) is the number of faces of $\Delta^m(B_{n-1})$ with 
$k-1$ elements, exactly one of which is a diameter, which also counts 
the $k$-element faces of 
$\Delta^m(D_n)$ with a pair of coinciding diameters.
Thus the first term of (\ref{piece2}) is the number of gray $k$-element faces 
with two or more diameters plus the number of dashed $k$-element faces 
with two or more diameters.
(This double-counts the faces with a pair of coinciding diameters.)

To complete the proof, it remains to show that the number of gray
(respectively dashed) 
$k$-element faces containing two or more diameters is 
\[
\binom{n-2}{k-2}\binom{(n-1)m+k}{k}.
\]
This is done by giving a bijection between gray faces with two or more diameters
and pairs 
\[
\left((\ep_1,\ep_2,\ldots,\ep_{n-2}),(a_1,a_2,\ldots,a_k)\right),
\]
where each $\ep_i$ is either 0 or 1, with $\ep_i=1$ exactly $k-2$ times, and 
the $a_j$ satisfy $1\le a_1\le a_2\le\cdots\le a_k\le(n-1)m+1$.
Given a gray face $F$ with $k$ vertices, including two or more diameters, 
the sequence of $a_j$'s is the set of indices of the starting vertices of 
the diameters or diagonal pairs in $F$, written in weakly increasing order.

For $k=2$, nothing more is needed to define the map, so assume $k>2$.
Now, take the smallest $j\in[k]$ such that $a_{j+1}-a_j>m$ (if such a $j$ exists),
where $a_{k+1}$ is understood to mean $a_1+(n-1)m+1$.
If the diagonal $[P_{a_j},P_{a_j+m+1}]$ is in 
$F$ then set $\ep_1=1$.  Otherwise set $\ep_1=0$.
The remaining $k-1$ of $k$ vertices of $F$ are a face in $\Delta^m(D_{n-1})$,
realized in the polygon $\Poly'$ obtained from $\Poly$ by deleting the vertices 
\[\pm P_{a_j+1},\pm P_{a_j+2},\ldots,\pm P_{a_j+m}.\]
By induction, determine the remaining $\ep_i$'s.

If there is no $j\in[k]$ with $a_{j+1}-a_j>m$, then we must have $k=n$, 
so set all the $\ep_i$'s to be 1.
Furthermore, in this case, all of the elements of $F$ must be diagonals.

To see that this is a bijection, we define an inverse.
Suppose we are given a pair 
$\left((\ep_1,\ep_2,\ldots,\ep_{n-2}),(a_1,a_2,\ldots,a_k)\right)$
and we wish to construct a face $F$.
If $k=2$, then $F$ consists of the gray diameter $[P_{a_1},-P_{a_1}]$
and the diameter $[P_{a_2},-P_{a_2}]$, the latter flavored in the unique
way which makes the two compatible.
Otherwise, take the smallest $j\in[k]$ such that $a_{j+1}-a_j>m$ (if such a $j$ exists).
If $\ep_1=1$ then 
use $\left((\ep_2,\ldots,\ep_{n-2}),(a_1,a_2,\ldots,\hat{a}_j,\ldots,a_k)\right)$
to inductively determine a $(k-1)$-element face $F'$ of $\Delta_m(D_{n-1})$
realized in the polygon $\Poly'$ described above.  Then $F$ is the 
face of $\Delta_m(D_n)$ obtained by adjoining the diagonal $[P_{a_j},P_{a_j+m+1}]$
and its symmetric diagonal to~$F'$.
If $\ep_1=0$, then use $\left((\ep_2,\ldots,\ep_{n-2}),(a_1,a_2,\ldots,a_k)\right)$
to inductively determine a $k$-element face $F'$ of $\Delta_m(D_{n-1})$
realized in~$\Poly'$.  Then $F$ is the face of $\Delta_m(D_n)$ which agrees
with~$F'$.

If no such $j$ exists, then $n$ must equal $k$.  Consider the collection
of $n$ diagonals whose starting vertices are given by the $a_i$'s.  By
Lemma~\ref{lem:all diams}, there is a unique way of assigning flavors to
each of the diameters in the collection so as to obtain a gray face of
$\Delta_m(D_n)$.
\end{proof}


\section{$h$-vectors}
\label{sec:h-vectors}

We define the \emph{$f$-polynomial} $F(\Phi,m,x)$
in the formal variable~$x$ by 
\[
F(\Phi,m,x)=\sum_{k=0}^n f_{n-k}(\Phi,m)\, x^k. 
\]
We then define the \emph{$h$-polynomial} $H(\Phi,m,x)$ by 
\[
H(\Phi,m,x)=F(\Phi,m,x-1).  
\]
The corresponding \emph{$h$-vector}
$(h_0,\dots,h_n)=(h_0(\Phi,m),\dots,h_n(\Phi,m))$ is then defined by 
\[
H(\Phi,m,x)=\sum_{k=0}^n h_{n-k}(\Phi,m)\, x^{k}. 
\]
For $m=1$, we recover the $h$-vector of the 
associahedron of type~$\Phi$, whose coefficients are the
(generalized) \emph{Narayana numbers} of type~$\Phi$. 
These Narayana numbers have been shown \cite{ga, panyushev,
  reiner-welker, sommers} 
to count various combinatorial objects associated with the root
system~$\Phi$, the Coxeter group~$W$, the corresponding semisimple Lie
algebra, etc. 
See \cite[Section~5.2]{pcmi} or \cite{Ath-TAMS} for a survey of these
results.

For a general~$m$, 
the numbers $h_k(\Phi,m)$ turn out to coincide with the $m$-analogues
of the Narayana numbers introduced and studied by
C.~A.~Athanasiadis~\cite{Ath-TAMS}. 
The latter numbers, herein denoted $H_k(\Phi,m)$, count various 
objects associated with a 
root system~$\Phi$
and a positive integer~$m$ according to certain statistics
whose values are recorded by~$k$.
We next sketch the definition of the numbers $H_k(\Phi,m)$, referring 
the reader to~\cite{Ath-TAMS} for
further details and references. 

The \emph{extended Catalan arrangement} 
$\operatorname{Cat}_\Phi^m$ (see \cite{Ath-Tokyo, postnikov-stanley})
is the collection of 
affine hyperplanes $H_{\alpha,j}$ defined by the equations $\br{\alpha,x}=j$
for $\alpha\in\Phi$ and $j=0,1,\ldots,m$.
The \emph{discrete torus} $T^m$ is the coroot lattice $\check{Q}$ of $\Phi$
modulo its dilation $(mh+1)\check{Q}$.
The group $W$ naturally acts on~$T^m$; the \emph{rank} of a $W$-orbit 
is the rank (as a reflection group) of the stabilizer of any point in the 
orbit.
The {\em root poset} of $\Phi$ is the partial order
on 
$\Phi_{>0}$ such that $\beta\le\gamma$ if and only if
$\gamma-\beta$ is a nonnegative linear combination of simple roots. 

The number $H_k(\Phi,m)$ counts:
\begin{itemize}
\item 
the regions $R$ of
$\operatorname{Cat}_\Phi^m$ in the fundamental
chamber such that $n-k$ walls 
of $R$ of the form $H_{\alpha,m}$ separate $R$
from the fundamental alcove; 
\item 
the $W$-orbits of rank $n-k$ in $T^m$;
\item 
collections of nested order filters in the root poset of~$\Phi$ 
satisfying certain technical conditions (see
\cite[(1.2)--(1.3)]{Ath-TAMS}) 
and having $n-k$ indecomposable elements of rank~$m$.
\end{itemize}

\begin{conjecture}\
\label{conj:we=athan}
For any finite crystallographic root 
system~$\Phi$, we have $h_k(\Phi,m)\!=\!H_k(\Phi,m)$. 
\end{conjecture}

Conjecture~\ref{conj:we=athan} can be verified for the 
classical series $ABCD$ by comparing the formulas
\eqref{eq:h-An}--\eqref{eq:h-Dn} below with their counterparts
in~\cite{Ath-TAMS}. 
This observation has already been made in~\cite{Tzanaki-1}. 
The conjecture has not been completely verified because the numbers $H_k(\Phi,m)$ have not been computed for the exceptional types.

The $h$-numbers $h_k(\Phi,m)$ are given by multiplicative formulas
very similar to the corresponding formulas for the face numbers. 

\begin{theorem}
\label{th:product-h}
For an irreducible root system $\Phi$ and an integer $0\le k\le n$, we
have 
\[
h_k(\Phi,m) 
= c_h(\Phi,k,m) \binom{n}{k} \prod_{\operatorname{level}(e_i)\leq k}
\frac{mh-e_i+1}{e_i+1} 
 \,, 
\]
where the factor
$c_h(\Phi,k,m)$ is a polynomial in $m$ given as follows:
\begin{itemize}
\item
in types $A_n$, $B_n$, $H_3$, and~$I_2(m)$, we have 
$c_h(\Phi,k,m)=1$;
\item
in type~$D_n$, we have 
\[
\qquad c_h(D_n,k,m)=
\begin{cases}
1 & \text{\ if $k\in\{0,1,n\!-\!1,n\}$;}\\
\displaystyle\frac{(n^2\!-\!n\!+\!k^2\!-\!k)(mn\!-\!m\!+\!1)}{kn(n\!-\!1)}\!-\!1
& \ \text{otherwise};
\end{cases}
\]
\item
in types $E_6$, $E_7$, $E_8$, $F_4$, and $H_4$, the factors $c_h(\Phi,k,m)$
are shown in Table~\ref{table:correction_h}. 
\end{itemize}
In particular, for the types $A_n$, $B_n$, and $D_n$, we have: 
\begin{align}
\label{eq:h-An}
h_k(A_n,m)&=\frac{1}{k+1}\binom{n}{k}\binom{(n+1)m}{k}\\
\label{eq:h-Bn}
h_k(B_n,m)&=\binom{n}{k}\binom{nm}{k}\\
\label{eq:h-Dn}
h_k(D_n,m)&=\binom{n}{k}\binom{(n-1)m}{k}+\binom{n-2}{k-2}\binom{(n-1)m+1}{k}.
\end{align}
\end{theorem}

\begin{table}[ht]
\begin{center}
\begin{tabular}{|c||c|c|c|c|c|c|}
\hline
$k$ & $D_8$ & $E_6$ & $E_7$   & $E_8$   & $F_4$ & $H_4$\\ \hline
&&&&&&\\[-3.5mm]
\hline
0 & $1$&$1$&$1$&$1$&$1$&$1$\\
\hline
1 & $1$&$1$&$1$&$1$&$1$&$1$\\
\hline
&&&&&&\\[\negtableskip]
2 & $\frac{203m-27}{56}$&$\frac{42m-8}{15}$&$\frac{63m-11}{18}$&$\frac{105m-17}{24}$&$\frac{78m-23}{36}$&$\frac{465m-149}{180}$\\[\tableskip]
\hline
&&&&&&\\[\negtableskip]
3 & $\frac{217m-53}{84}$&$\frac{18m-5}{8}$&$\frac{27m-7}{10}$&$\frac{45m-11}{14}$&$1$&$1$\\[\tableskip]
\hline
&&&&&&\\[\negtableskip]
4 & $\frac{119m-39}{56}$&$\frac{30m-13}{18}$&$\frac{63m-23}{30}$&$\frac{10350
  m^2-6675 m+1084}{1575}$&$1$&$1$\\[\tableskip]
\hline
&&&&&&\\[\negtableskip]
5 & $\frac{133m-51}{70}$&$1$&$\frac{207m-103}{126}$&$\frac{30m-13}{15}$&&\\[\tableskip]
\hline
&&&&&&\\[\negtableskip]
6 & $\frac{301m-125}{168}$&$1$&$1$&$\frac{195m-107}{120}$&&\\[\tableskip]
\hline
7 & $1$&&$1$&$1$&&\\
\hline
8 & $1$&&&$1$&&\\
\hline
\end{tabular}
\end{center}
\caption{The factors $c_h(\Phi,k,m)$}
\label{table:correction_h}
\end{table}

\begin{proof}
As in the case of Theorem~\ref{th:product-f}, the type $I_2(a)$ is
an easy calculation, and the types $E_6$, $E_7$, $F_4$, $G_2$, $H_3$, and $H_4$ can
be verified by computer. 
The formulas in the classical types $ABD$ can be derived from their counterparts in
Theorem~\ref{th:product-f} using Chu-Vandermonde summation. 
Details are \hbox{omitted.} 
\end{proof}

\section{Euler characteristic}
\label{sec:euler}

Recall that the \emph{reduced Euler characteristic} of a simplicial
complex~$\Delta$ that has $f_k$ faces of dimension $k-1$, for
$k=0,1,\dots,n$ (including one empty face for $f_0=1$), is defined by 
\[
\tilde\chi(\Delta) = \sum_{k\geq 0} (-1)^{k-1} f_k\,. 
\]
\pagebreak[3]
Thus $\tilde\chi(\Delta)=\chi(\Delta)-1$, where $\chi(\Delta)$ is the
usual Euler characteristic. 
It is well-known that the Euler characteristic of an $(n-1)$-dimensional
complex is $(-1)^{n-1}h_n$.
Using the $k=n$ case of Theorem~\ref{th:product-h} and the fact that
$e_i=h-e_{n-i}$, we obtain: 

\begin{proposition}
\label{prop:euler-recip}
The reduced Euler characteristic of $\Delta^m(\Phi)$ is given by 
\begin{equation}
\label{eq:euler-char-phi}
\tilde\chi(\Delta^m(\Phi)) = (-1)^{n-1} N(\Phi,m\!-\!1). 
\end{equation}
Thus, it is equal, up to a sign, to the number of maximal simplices 
in the complex~$\Delta^{m-1}(\Phi)$. 
\end{proposition}

Formula~\eqref{eq:euler-char-phi} can be restated as 
\begin{equation}
\label{eq:euler-char-phi-restated}
\sum_{k=0}^n  (-1)^{n-k} f_k(\Phi,m) = N(\Phi,m\!-\!1). 
\end{equation}

In the special case $m=1$, 
Proposition~\ref{prop:euler-recip} is a corollary of the following
statement proved in~\cite{ga}: 
the simplicial complex $\Delta^1(\Phi)$ is homeomorphic to an
$(n-1)$-dimensional sphere. 
Note that $N(\Phi,0)=1$ (see Example~\ref{example:m=0}). 

\begin{remark}
Our proof of Proposition~\ref{prop:euler-recip} relies on a
type-by-type calculation. 
Is there a conceptual and/or geometric proof? 
For example, there may be a way to remove
$N(\Phi,m\!-\!1)$ open facets from~$\Delta^m(\Phi)$ and obtain a
contractible complex.
\end{remark}

The following conjecture was suggested by our discussions with
Hugh Thomas in December 2004. 

\begin{conjecture} 
\label{conj:thomas}
The complex $\Delta^m(\Phi)$ is shellable, hence Cohen-Macaulay, and
homotopy equivalent to a wedge of spheres. 
\end{conjecture} 

Conjecture~\ref{conj:thomas} has been proved by E.~Tzanaki
\cite[Theorem~1.4]{Tzanaki} in types $A_n$ and~$B_n\,$. 

\section{Reciprocal face numbers}
\label{sec:recip}

Inspired by the construction of ``positive clusters'' in~\cite{ga}, 
we define the \emph{reciprocal face numbers} $f_k^+(\Phi,m)$ as follows.
Since $f_k(\Phi,m)$ is a polynomial in~$m$, 
one can extend the definition of $f_k(\Phi,m)$ to the
non-positive integer values of~$m$. 
Let us define
\begin{equation}
\label{eq:def-fk+}
f_k^+(\Phi,m) = (-1)^k f_k(\Phi,-m\!-\!1) .
\end{equation}
For $k=n$, we set
\begin{equation}
\label{eq:def-N+}
N^+(\Phi,m)=f_n^+(\Phi,m) = (-1)^n N(\Phi,-m-1) .
\end{equation}
Our choice of notation~$N^+$ is explained by
Proposition~\ref{prop:incl-excl-parabolic} below. 

The definition \eqref{eq:def-N+} combined with 
Propositions~\ref{prop:fk-recur} and~\ref{pr:product-N} 
yields the following formulas. 

\begin{corollary}
\label{cor:recurr-N+}
For a reducible root system $\Phi=\Phi_1\times \Phi_2$,
\begin{equation}
\label{eq:N+-reducible}
N^+(\Phi_1\times \Phi_2,m)=N^+(\Phi_1,m)\,N^+(\Phi_2,m).
\end{equation}
If $\Phi$ is irreducible, then
\begin{equation}
\label{eq:recurr-N+}
N^+(\Phi,m)=\frac{(m\!+\!1)h\!-\!2}{2n}\sum_{i\in I}N^+(\Phi_{\br{i}},m).
\end{equation}
\end{corollary}

\begin{corollary}
\label{cor:product-N+}
For an irreducible root system $\Phi$,
\begin{equation}
\label{eq:product-N+}
N^+(\Phi,m)=\prod_{i}\frac{mh+e_i-1}{e_i+1} 
\end{equation}
\end{corollary}

We note that the recurrence \eqref{eq:recurr-N+}
is new even for $m=1$, with the numbers $N^+(\Phi)=N^+(\Phi,1)$ counting
positive clusters in~$\Delta(\Phi)$ (cf.\ \cite[Proposition~3.9]{ga}). 

Analogues of \eqref{eq:recurr-N+} and \eqref{eq:product-N+} 
for the general numbers $f_k^+(\Phi,m)$
can be obtained   directly from \eqref{eq:def-fk+}
and similar formulas for the face numbers $f_k(\Phi,m)$. 

\begin{proposition}
\label{prop:incl-excl}
We have 
\begin{equation}
\label{eq:N:incl-excl}
N(\Phi,m)=\sum_{J\subseteq I} N^+(\Phi_J,m),
\end{equation}
where, by convention, $N^+(\Phi_J,m)=1$ for the empty subset
$J=\emptyset$. 
\end{proposition} 

\begin{proof}
As shown in \cite[Corollary~1.3]{Ath} (cf.\ discussion in
Section~\ref{sec:h-vectors}), 
the number $N(\Phi,m)$ (resp., $N^+(\Phi,m)$) counts all (resp., bounded)
regions of the extended Catalan arrangement
$\operatorname{Cat}_\Phi^m$ contained in the fundamental chamber
of the corresponding Coxeter arrangement. 
Then \eqref{eq:N:incl-excl} follows, for example, from
\cite[Lemma~5.3]{Ath-Tz}. 
\end{proof}

The following generalization of \cite[Proposition~3.9]{ga} 
provides a direct combinatorial interpretation of the numbers
$N^+(\Phi,m)$ in terms of generalized cluster complexes. 

\begin{proposition}
\label{prop:incl-excl-parabolic}
The number of maximal simplices in $\Delta^m(\Phi)$ which only 
involve positive colored roots is equal to~$N^+(\Phi,m)$. 
\end{proposition} 

\begin{proof}
Let $\tilde N^+(\Phi,m)$ denote the number of maximal simplices in
$\Delta^m(\Phi)$ consisting of positive colored roots. 
Theorem~\ref{th:restrict} implies that 
$N(\Phi,m)=\sum_{J\subseteq I} \tilde N^+(\Phi_J,m)$,  
and the claim  follows by Proposition~\ref{prop:incl-excl}. 
\end{proof}


\begin{remark}
It is natural to suggest that the relationship \eqref{eq:def-N+}
between the quantities $N(\Phi,m)$ and $N^+(\Phi,m)$ is an instance of
Ehrhart reciprocity~\cite[Theorem~4.6.26]{EC1}. 
This suggestion has been confirmed by C.~A.~Athanasiadis and
E.~Tza\-naki in \cite[Section~7]{Ath-Tz}. 
\end{remark}

\section{Combinatorics of Coxeter-theoretic invariants}
\label{sec:coxeter-inv}

A \emph{Coxeter diagram} is an undirected graph~$G$ 
without loops and multiple edges, 
in which every edge is labeled by an integer~$\geq 3$. 
(In general Coxeter group theory, labels equal to~$\infty$ are also
allowed, but we restrict our attention to integer labels.) 
By convention, each missing edge (i.e., each pair of vertices in~$G$
not connected by an edge) is thought of as an edge labeled by~$2$. 
The number of vertices in a Coxeter diagram is called its
\emph{rank}. 


Let $G$ be a Coxeter diagram of rank $n$ for an irreducible finite
Coxeter group~$W$. 
There is a host of Coxeter-theoretic invariants associated
with~$G$ and the corresponding root system~$\Phi$:
\begin{itemize}
\item
the exponents $e_1,\dots,e_n$; 
\item
the Coxeter number $h=\frac{2}{n}\displaystyle\sum e_i$;
\item
the number of roots $|\Phi|=nh$; 
\item
the order of the group $|W|=\displaystyle\prod (e_i+1)$; 
\item
the generalized Catalan number $N(\Phi)=\displaystyle\prod
\frac{h+e_i+1}{e_i+1}$; 
\item
the face numbers $f_k(G,m)$ of the simplicial
complex~$\Delta^m(\Phi,m)$; 
\end{itemize}
and many others. 
Our goal in this section is to develop purely combinatorial procedures
for computing all these invariants directly from~$G$,
i.e., procedures that do not explicitly involve the root
system~$\Phi$, the Coxeter group~$W$, or any other algebraic or
Lie-theoretic concepts. 

\begin{remark}
Since all quantities listed above can be expressed
in terms of the exponents $e_1,\dots,e_n$, 
computing the exponents would accomplish the task. 
However, determining them directly
from the Coxeter diagram~$G$ is not at all straightforward. 

In turn, the computation of exponents can be reduced to the
calculation of a Coxeter number, as follows. 
Use Proposition~\ref{prop:fk-recur}
to recursively compute the generalized Fuss number 
\[
N(G,m)=N(\Phi,m)=\prod \frac{mh+e_i+1}{e_i+1}.  
\]
Then, viewing $N(G,m)$ as a polynomial in~$m$, 
find its roots $-\frac{e_i+1}{h}$ and determine the
exponents~$e_i$. 
\end{remark}

\begin{remark}
As a rule of thumb, any use of a ``Coxeter number'' or ``exponents'' in the
reducible case should be viewed with skepticism.  In particular, the
formulas above in this section expressing other invariants in
terms of the Coxeter number and/or the exponents make no sense in the
reducible case. In what follows, we do not consider Coxeter numbers or
exponents for any disconnected diagrams. 
\end{remark}

We will describe several combinatorial algorithms for computing
invariants of Coxeter diagrams associated with finite Coxeter groups. 
Each of these algorithms can then be used to define ``fake''
Coxeter-theoretic invariants of more general Coxeter diagrams.
Such generalizations are discussed in Section~\ref{sec:fake}.

All our algorithms are recursive:
the computation of invariants of~$G$ relies upon prior
calculation of similar invariants for induced subgraphs of~$G$. 
The base of recursion is provided by the special
cases~$n\leq 2$, together with setting   
$f_0(G,m)=1$ for any~$G$. 

For $n=0$, $G=\emptyset$, we have 
\[
N(G,m)=N^+(G,m)=1. 
\]

For $n=1$, $G=(\bullet)$ (see Example~\ref{example:f:n=1}),  we
postulate: 
\begin{align}
\nonumber
h&=2,\\
\nonumber
e_1&=1,\\
\label{eq:f1:n=1}
N(G,m)=f_1(G,m) &= m+1, \\
\label{eq:N+:n=1}
N^+(G,m)&=m. 
\end{align}

For $n=2$, $G=(\bullet\overunder{a}{}\bullet)$ (see Example~\ref{example:f:n=2}),  we
postulate: 
\begin{align}
\nonumber
h&=a,\\
\nonumber
\{e_1,e_2\}&=\{1,a-1\},\\
\label{eq:f1:n=2}
f_1(G,m) &= am+2, \\
\label{eq:f2:n=2}
N(G,m)=f_2(G,m) &= \frac{(am+2)(m+1)}{2},\\
\label{eq:N+:n=2}
N^+(G,m)&=\frac{(am+a-2)m}{2}.
\end{align}
We note that formulas \eqref{eq:f1:n=2}--\eqref{eq:N+:n=2} 
hold just as well in the reducible
case~$a=2$, where $G$ is of type $A_1\times A_1$,
and we have $f_1(A_1\times A_1)=2m+2$ and $N(A_1\times A_1)=(m+1)^2$.

We next present four different combinatorial procedures, each
of which determines the invariants of a connected
Coxeter diagram $G$ of rank $n\geq 3$ assuming that such
invariants for all proper induced connected subgraphs of~$G$ are
already known. 

\addtocontents{toc}{\SkipTocEntry}
\subsection{The Euler characteristic method} 
\label{sec:euler-char-method}

We start with the basic recurrence \eqref{eq:recur-fk}: 
\begin{equation}
\label{eq:diffeq-f}
f_k(G,m) = \frac{mh+2}{2k} \sum_{G'\lessdot G} f_{k-1}(G',m)\,,
\end{equation}
where $k\in\{1,\dots,n\}$ and the summation is over all
induced subgraphs $G'\subset G$ obtained by removing one vertex
from~$G$. 
Whenever $G'$ is disconnected,  with connected components $G_1,\dots,G_s$, 
we use the formula 
\begin{equation}
\label{eq:fk-disconnected}
f_{k-1}(G',m)=\sum_{k_1+\cdots+k_s=k-1} \ \prod_{i=1}^s f_{k_i}(G_i,m) 
\end{equation}
(see  \eqref{eq:fk-reducible}) to calculate $f_{k-1}(G',m)$. 

We note that the number~$h=h(G)$ is undetermined as of yet. 
To find the value of~$h$, 
we use~\eqref{eq:euler-char-phi-restated}: 
\begin{equation}
\label{eq:euler-char-phi-G}
\sum_{k=0}^n  (-1)^{n-k} f_k(G,m) = f_n(G,m\!-\!1). 
\end{equation}
Substituting \eqref{eq:diffeq-f} into \eqref{eq:euler-char-phi-G}
gives a linear equation for~$h$:
\begin{align}
\label{eq:euler-char-phi-G-recur}
&(-1)^n\!+\sum_{G'\lessdot G} \Bigl(
-\frac{(m\!-\!1)h\!+\!2}{2n}f_{n-1}(G',m\!-\!1)
\!+\! (mh\!+\!2)\sum_{k=1}^{n}  \frac{(-1)^{n-k}}{2k} f_{k-1}(G',m) 
\Bigr)
= 0. 
\end{align}
\emph{A~priori}, solving this equation for~$h$ 
results in a rational function in~$m$.
However, if $G$ is a Coxeter diagram for a finite Coxeter group
(and in some other cases as well), the answer miraculously turns out
to be a constant---in fact, a positive integer.

\begin{example}[$n=3$] 
\label{example:euler-n=3}
Let
$G$ be a connected Coxeter diagram of 
rank~$3$ with edge labels $a_1,a_2,a_3\geq 2$, as 
shown in Figure~\ref{fig:tri_diag}.

\begin{figure}[ht]
\centerline{\scalebox{1}{\epsfbox{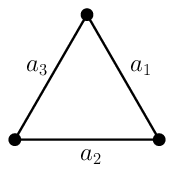}}}
\caption{}
\label{fig:tri_diag}
\end{figure}

We denote $a=a_1+a_2+a_3\,$.
The values of $a$ corresponding to finite irreducible Coxeter groups
are shown in Table~\ref{table:a-for-A3B3H3}. 

\begin{table}[ht]
\begin{center}
\begin{tabular}{|c|c|c|c|}
\hline
type of $G$ & $A_3$ & $B_3$ & $H_3$ \\\hline
&&&\\[-3.5mm]
\hline
$a$& $8$ & $9$ & $10$ \\\hline
\end{tabular}
\end{center}
\vspace{.2in}
\caption{}
\label{table:a-for-A3B3H3}
\end{table}

The ``maximal parabolic'' subgraphs $G'\lessdot  G$ are
of types $I_2(a_1)$, $I_2(a_2)$, and $I_2(a_3)$, respectively.
Their contributions to the outer sum in
~\eqref{eq:euler-char-phi-G-recur} are of the form 
\begin{align*}
&-\frac{(m\!-\!1)h+2}{6} f_2(I_2(a_i),m-1)
+(mh+2)\sum_{k=1}^3  \frac{(-1)^{3-k}}{2k}
f_{k-1}(I_2(a_i),m)\\
=&-\frac{(m\!-\!1)h+2}{6} \frac{(a_i(m-1)+2) m}{2}\\
&\qquad\qquad +(mh+2)\left(
\frac{1}{2}
-\frac{1}{4}(a_i m+2)
+\frac{1}{6}\frac{(a_i m+2)(m+1)}{2}\right)\\
=&\,\frac{-a_i mh-2a_i m+4mh+4}{12}\,.
\end{align*}
Thus \eqref{eq:euler-char-phi-G-recur} becomes
\[
-1 
-\frac{(a_1+a_2+a_3)m(h+2)}{12} +mh+1
=0
\]
yielding
\begin{equation}
\label{eq:h:n=3}
h=\frac{2a}{12-a}\,.   
\end{equation}
We can now use this value of $h$ to compute the face numbers
by means of~\eqref{eq:diffeq-f}:
\begin{align}
f_1(G,m) &= \frac{3(am-a+12)}{12-a}, \\
f_2(G,m) &= \frac{(am-a+12)(am+6)}{2\,(12-a)}, \\
\label{eq:N(G),n=3}
N(G,m)=f_3(G,m) &= \frac{(m+1)(am+6)(am-a+12)}{6\,(12-a)}
\end{align}
The roots of $N(G,m)$ are 
$\{\mu_1,\mu_2,\mu_3\}=\{-1, -\frac{6}{a},-\frac{12-a}{a}\}$.  
Since 
$\mu_i=-\frac{e_i+1}{h}$, we recover the exponents by 
$e_i=-h\mu_i-1$: 
\begin{equation}
\label{eq:e1e2e3}
\{e_1,e_2,e_3\}
=\Bigl\{1, \frac{a}{12-a}, \frac{3a-12}{12-a}\Bigr\}. 
\end{equation}
\end{example}
We also use \eqref{eq:def-N+} to get
\begin{equation}
\label{eq:N+:n=3}
N^+(G,m)=\frac{m(am+a-6)(am+2a-12)}{6(12-a)}. 
\end{equation}

\addtocontents{toc}{\SkipTocEntry}
\subsection{Symmetry-based method}

This method does not require a recursive computation of all 
polynomials~$f_k(G,m)$, but only~$N(G,m)$. 
It exploits the fact that the set of exponents is invariant 
with respect to the reflection
\[
e_i \mapsto h-e_i \,. 
\]
Equivalently, the set $\{\mu_i\}$ of (simple) roots of the polynomial
$N(G,m)$ is invariant 
with respect to the reflection
\begin{equation}
\label{eq:mu-reflection}
\mu_i \mapsto -\frac{h+2}{h}-\mu_i \,. 
\end{equation} 

Recall once again the recurrence
\begin{equation}
\label{eq:diffeq-N}
N(G,m) = \frac{mh+2}{2n} \sum_{G'\lessdot G} N(G',m)\,. 
\end{equation}
It follows by induction on~$n$ that $N(G,-1)=0$ for any
non-empty~$G$. 
Hence 
\begin{equation}
\label{eq:recur-Q}
  Q(G,m)= \frac{1}{m+1}\sum_{G'\lessdot G} N(G',m)
\end{equation}
is a polynomial in~$m$ (of degree~$n-2$). Furthermore,
its set of roots is still invariant under the
reflection~\eqref{eq:mu-reflection},
since we removed two symmetric roots $-1$ and $-\frac2h$
from~$N(G,m)$. 

Note that the computation of $Q(G,m)$ by the formula \eqref{eq:recur-Q} 
does not require the value of~$h$. So we find the polynomial~$Q(G,m)$, 
and then determine the average of its roots.
This average is equal to~$-\frac{h+2}{2h}$, 
so the negative sum of the $n-2$ roots of $Q(G,m)$ is equal to
\begin{equation}
\label{eq:coeff/coeff}
  \frac{\text{coefficient of $m^{n-3}$ in~$Q(G,m)$}}{\text{coefficient of
      $m^{n-2}$ in~$Q(G,m)$}}
= \frac{(n-2)(h+2)}{2h}.
\end{equation}
We then find $h$ from this linear equation, and determine the
exponents. 

\begin{example}[$n=3$] 
\label{example:symmetry-n=3}
Let $G$ be a connected Coxeter diagram of 
rank~$3$ with edge labels $a_1,a_2,a_3$ (cf.\
Example~\ref{example:euler-n=3}). 
As before, we use the notation $a=a_1+a_2+a_3\,$. 
Combining \eqref{eq:recur-Q} with \eqref{eq:f2:n=2}, we obtain:
\[
Q(G,m)= \frac{1}{m+1}\sum_{i=1}^3 \frac{(a_im+2)(m+1)}{2}
=\frac{am+6}{2}. 
\]
Hence \eqref{eq:coeff/coeff} becomes
\[
\frac{6}{a}=\frac{h+2}{2h}, 
\]
implying \eqref{eq:h:n=3}. 
We conclude that the roots of $N(G,m)$ are $-\frac{6}{a}$ 
(the only root of~$Q(G,m)$), 
$-1$ (always a root), and $-\frac{12-a}{a}$ 
(obtained from~$-1$ using~\eqref{eq:mu-reflection}).
The exponents \eqref{eq:e1e2e3} are then found as in
Example~\ref{example:euler-n=3}. 
\end{example}

\addtocontents{toc}{\SkipTocEntry}
\subsection{Reciprocity-based method}
\label{sec:reciprocity-based} 

This method uses the inclusion-exclusion formula for the reciprocal
face numbers $N^+(G,m)$ as the additional equation needed to
recursively compute the numbers~$N(G,m)$ if the Coxeter number~$h$ is
unknown. For the simplest version of this method, 
one only needs to recursively compute the
numbers $N(G)=N(G,1)$ and $N^+(G)=N^+(G,1)=(-1)^n N(G,-2)$.
That is, this version only involves the numbers $f_k(G,m)$ for $k=n$
and $m\in\{1,-2\}$. 

The idea is very simple. Proposition~\ref{prop:incl-excl}
and formula \eqref{eq:diffeq-f} give: 
\begin{align}
\label{eq:recip-method-1}
N(G) &= {\frac{h+2}{2n}} \sum_{G'\lessdot G} N(G'),\\
\label{eq:recip-method-2}
N^+(G) &= {\frac{h-1}{n}} \sum_{G'\lessdot G} N^+(G'), \\
\label{eq:recip-method-3}
N(G) &= \sum_{H\subset G} N^+(H) ,
\end{align}
where $H$ ranges over all induced subgraphs of~$G$ (including~$G$ and
the empty subgraph). 
Since the values of $N(H)$ and $N^+(H)$ for all proper subgraphs
of~$G$ are presumed known, it remains to find the three unknowns
$N(G)$, $N^+(G)$, and~$h$ from the three linear equations
\eqref{eq:recip-method-1}--\eqref{eq:recip-method-3}.

\begin{example}[$n=3$] 
\label{example:recip-n=3}
Let $G$ be a connected Coxeter diagram of 
rank~$3$ with edge labels $a_1,a_2,a_3$ (cf.\
Examples~\ref{example:euler-n=3} and~\ref{example:symmetry-n=3}). 

The values of $N(H)$ and $N^+(H)$ for 
graphs $H$ of rank $1$ or $2$ 
are given by \eqref{eq:f1:n=1} and~\eqref{eq:f2:n=2}
 with $m=1$ and $m=-2$: 
\[
\begin{array}{lll}
H=(\bullet) & \hspace{.5in} N(H)=2 & \hspace{.5in} N^+(H)=1 \\
H=(\bullet\overunder{b}{}\bullet) & \hspace{.5in} N(H)=b+2 & \hspace{.5in} N^+(H)=b-1
\end{array}
\]
We continue to use the notation $a=a_1+a_2+a_3\,$. 
Equations \eqref{eq:recip-method-1}--\eqref{eq:recip-method-3} become:
\begin{align*}
N(G) &= {\frac{h+2}{6}} \,(a+6),\\
N^+(G) &= {\frac{h-1}{3}} \,(a-3), \\
N(G) &= N^+(G)+(a-3)+3+1 .
\end{align*}
Solving these equations for $h$, $N(G)$, and $N^+(G)$ yields 
$h=\frac{2a}{12-a}$, $N(G)=\frac{4(a+6)}{12-a}$,
and $N^+(G)=\frac{(a-3)(a-4)}{12-a}$,
in agreement with \eqref{eq:h:n=3}, \eqref{eq:N(G),n=3}, and
\eqref{eq:N+:n=3} (with $m=1$). 
The exponents are found as before. 
\end{example}

The more general version of the reciprocity-based method retains the
parameter~$m$. 
Let $G$ be a connected Coxeter diagram of rank~$n\geq 3$. 
As before, we start with the relevant versions of the basic
recurrence~\eqref{eq:diffeq-f}, namely \eqref{eq:diffeq-N} and 
\begin{align}
\label{eq:recip-method-5}
N^+(G,m) &= {\frac{h(m+1)-2}{2n}} \sum_{G'\lessdot G} N^+(G',m).
\end{align}
Substituting (\ref{eq:recip-method-3}) into both sides of 
 \eqref{eq:diffeq-N}, we obtain:  
\begin{align*}
\sum_{H\subset G} N^+(H,m) &= \frac{mh+2}{2n}\sum_{H\subset G}
(n-|H|)N^+(H,m),  
\end{align*}
where $|H|$ denotes the rank of~$H$. 
Solving for~$N^+(G,m)$, we get 
\begin{align}
\label{eq:recip-method-6}
N^+(G,m) & =\sum_{H\subsetneq G} \frac{mh}{2n}\, (n-|H|)\, N^+(H,m).
\end{align}
Equating the right-hand sides of 
(\ref{eq:recip-method-6}) and (\ref{eq:recip-method-5}) 
and solving for~$h$ yields
\begin{align}
\label{eq:recip-method-7}
h & =2\frac{(n-2)\sum_{G'\lessdot G}N^+(G',m)+\sum_{H\ll G} |H|N^+(H,m)}
        {m\sum_{H\ll G}(n-|H|) N^+(H,m)-\sum_{G'\lessdot G}N^+(G',m)},
\end{align}
where  $H\ll G$ means that $H\subset G$ and $|H|<|G|-1$. 

\addtocontents{toc}{\SkipTocEntry}
\subsection{The invariants $M(G)$} 
\label{sec:M(G)}

We next present a simplified version of the reciprocity-based method
that exploits the fact that the Coxeter diagrams of finite Coxeter groups
are \emph{sparse}, i.e., they have far fewer edges than~$\binom{n}{2}$. 
The method is based on recursive computation of the following invariant. 

\begin{definition}
\label{def:M(G)} 
For a non-empty Coxeter graph~$G$ of rank~$n$, define
\begin{equation} 
\label{eq:def-M(G)}
M(G)= n \lim_{m\to 0} \frac{N^+(G,m)}{m}. 
\end{equation} 
If $G$ is empty, set $M(G)=0$. 
\end{definition}

\begin{remark}
\label{rem:M(G)-disconn}
Note that \eqref{eq:N+:n=1}, \eqref{eq:N+:n=2}, and
\eqref{eq:recip-method-5} imply $N^+(G,0)=0$, so $M(G)$ is well
defined. 
For the same reasons, $M(G)=0$ if $G$ is disconnected.
\end{remark}

\pagebreak[3]

\begin{corollary}
\label{cor:chapoton-plein}
If $G$ is the Coxeter diagram of an irreducible finite Coxeter group
with exponents $e_1,\dots,e_n$ and Coxeter number~$h$, then 
\begin{equation}
\label{eq:chapoton-plein}
M(G)=\frac{nh}{2}\prod_{e_i\neq 1}\frac{e_i-1}{e_i+1} .
\end{equation}
\end{corollary}

\begin{proof}
Follows from \eqref{eq:product-N+} and~\eqref{eq:def-M(G)}. 
\end{proof}

\begin{table}[ht]
\begin{center}
\begin{tabular}{|c||c|c|c|c|c|c|c|c|c|c|}
\hline
type of $G$ & $A_n$  & $B_n$    & $D_n$    & $E_6$ & $E_7$   & $E_8$   & $F_4$ & $H_3$& $H_4$ & $I_2(a)$ \\ \hline
&&&&&&&&&&\\[-4.5mm]
\hline
$M(G)$ & $1$ & $n$ & $n-2$ & $7$ & $16$ & $44$ & $10$ & $8$ & $42$ & $a-2$ \\\hline
\end{tabular}

\bigskip

\end{center}
\caption{Values of $M(G)$ for finite Coxeter groups}
\label{table:M(G)}
\end{table}

The values of $M(G)$ for all finite irreducible Coxeter groups  are
listed in Table~\ref{table:M(G)}. 
All of these values are integers. 
An explanation of this integrality is provided by the following
interpretation of the invariants~$M(G)$, due to F.~Chapoton. 

\pagebreak[3]

\begin{proposition}[\cite{chapoton-plein}]
\label{prop:plein}

Let $W$ be a finite irreducible Coxeter group with Coxeter
diagram~$G$. 
Then the number of reflections in~$W$ 
which do not belong to any proper parabolic subgroup~$W_J$ 
is equal to~$M(G)$. 
Accordingly, the number of all reflections in~$W$ is equal to 
\begin{equation}
\label{eq:plein-incl-excl}
\frac{nh}{2}=\sum_{H\subset G} M(H)\,. 
\end{equation}
\end{proposition}

Proposition~\ref{prop:plein}  can be restated as saying that in a
finite root system~$\Phi$, the number of positive roots whose simple
root expansion has full support (i.e., involves every simple root) is equal
to~$M(G)$, where $G$ is the corresponding Coxeter diagram. 

\begin{proof}
The inclusion-exclusion relation~\eqref{eq:plein-incl-excl} uniquely
defines the numbers~$M(G)$, and is satisfied by the number of
reflections not lying in a proper parabolic subgroup.
Consequently, it suffices to prove~\eqref{eq:plein-incl-excl} in order
to prove the Proposition. 

Substituting \eqref{eq:N:incl-excl} into both sides of
\eqref{eq:recurr-N}, we obtain: 
\[
\sum_{H\subset G} N^+(H,m)
=\frac{mh+2}{2n} \sum_{H\subset G} (n-|H|)\, N^+(H,m)\,,
\]
implying
\[
\frac{mh}{2}\sum_{H\subset G} N^+(H,m)
=\frac{mh+2}{2n} \sum_{H\subset G} |H|\, N^+(H,m)\,. 
\]
Dividing by $m$ and substituting $m=0$, we get (the only surviving
term on the left-hand side is the one for $H=\emptyset$):
\[
\frac{h}{2}=\frac{1}{n}\sum_{H\subset G} M(H)\,,
\]
and \eqref{eq:plein-incl-excl} follows. 
\end{proof}

It would be interesting to find a direct proof of
Proposition~\ref{prop:plein} based on Definition~\ref{def:M(G)} 
and one of the combinatorial interpretations of $N^+(G,m)$,  
rather than on calculations involving product formulas 
and/or recurrence relations. 

It follows from \eqref{eq:recurr-N+} that the
numbers $M(G)$ satisfy the recurrence
\begin{equation}
\label{eq:recurr-M}
M(G)=\frac{h\!-\!2}{2(n-1)}\sum_{G'\lessdot G} M(G').
\end{equation}
There is also a more complicated recurrence for these numbers that
does not require knowing the Coxeter number~$h$. 


\begin{lemma}
\label{lem:M(G)-recur}
Let $G$ be the Coxeter diagram of a finite irreducible Coxeter group
of rank $n\geq 3$. 
Set 
\begin{equation}
\label{eq:M(G)-ingredients}
\Sigma_1=\sum_{G'\lessdot G} M(G'), \qquad
\Sigma_2=\sum\limits_{
\begin{array}{c}\ \\[-.25in] \scriptstyle H\subset G\\[-.05in] 
                \scriptstyle 2\leq |H|\leq n-1\end{array}} 
M(H)\,. 
\end{equation}
Then 
\begin{align}
\label{eq:M(G)-rec}
M(G)&=\frac{\Sigma_1\,\Sigma_2}{n(n-1)-\Sigma_1}\,. 
\end{align}
\end{lemma}

\begin{proof}
Identities \eqref{eq:plein-incl-excl} and \eqref{eq:recurr-M} 
can be written as
\begin{align}
\label{eq:h-via-M}
h&=\frac{2}{n} (M(G)+\Sigma_2+n)\,,\\
\label{eq:M-via-h}
M(G)&=\frac{h\!-\!2}{2(n-1)}\Sigma_1\,.
\end{align}
Substituting \eqref{eq:h-via-M} into \eqref{eq:M-via-h}
yields~\eqref{eq:M(G)-rec}. 
\end{proof}

Formulas \eqref{eq:M(G)-ingredients}--\eqref{eq:M(G)-rec} can be used
to recursively compute the values~$M(G)$.
In view of Remark~\ref{rem:M(G)-disconn},
the sums in \eqref{eq:M(G)-ingredients} 
can be restricted to \emph{connected} subgraphs $G'$ and~$H$,
simplifying the calculations. 

The Coxeter numbers are then recovered using
\eqref{eq:plein-incl-excl}
or, equivalently,~\eqref{eq:h-via-M}. 

\begin{example}
\label{example:K3-M(G)}
For $G$ of type $I_2(a)$, one has $M(G)=a-2$. 
For a $3$-vertex diagram with edge labels $a_1,a_2,a_3$  
(cf.\ Examples~\ref{example:euler-n=3}, \ref{example:symmetry-n=3} 
and~\ref{example:recip-n=3}), with $a=a_1+a_2+a_3$, 
we get $\Sigma_1=\Sigma_2=a-6$. 
Substituting this into~\eqref{eq:M(G)-rec}--\eqref{eq:h-via-M},
we obtain 
\begin{equation}
\label{eq:M(G),n=3}
M(G)=\frac{(a-6)^2}{12-a}\,, \qquad h(G)=\frac{2a}{12-a}\,, 
\end{equation}
matching the results of our earlier calculations. 
\end{example}

\begin{example}
\label{example:4-chain}
Now let 
\[
G = (\bullet\overunder{b_1}{}\bullet\overunder{b_2}{}\bullet\overunder{b_3}{}\bullet)
\]
\pagebreak[3]
be a linear Coxeter diagram of rank~$4$.
(Here $b_1,b_2,b_3$ are integers~$\geq 3$.) 
Using~\eqref{eq:M(G),n=3}, we compute: 
\begin{align*}
\Sigma_1&=
\frac{(b_1+b_2-4)^2}{10-b_1-b_2}+
\frac{(b_2+b_3-4)^2}{10-b_2-b_3}\,,\\[.1in] 
\Sigma_2&=b_1+b_2+b_3-6+\Sigma_1\,. 
\end{align*}
The values of $M(G)$ and $h(G)$ are then obtained using
\eqref{eq:M(G)-rec}--\eqref{eq:h-via-M}. 
The resulting formula for the Coxeter number $h(G)$
is a cumbersome rational expression in $b_1,b_2,b_3$.
A diligent reader might want to verify that,
in the special cases where $G$ is of types~$A_4$, $B_4$, $F_4$,
or~$H_4$, the outcomes coincide with the known values of Coxeter
numbers. 
\end{example}

\section{Fake Coxeter invariants}
\label{sec:fake}

Each of the methods described in Section~\ref{sec:coxeter-inv}
can, in principle, be applied to a very general Coxeter diagram~$G$. 
The resulting ``fake Coxeter invariants'' 
may or may not be meaningful.  
Be that as it may, it is tempting to calculate these invariants in the
simplest non-finite cases; 
to catalogue the infinite types where such invariants are ``nice;''
and possibly find intrinsic interpretations for them.
One also wonders what the limits of applicability are 
for each of the methods of Section~\ref{sec:coxeter-inv}, 
and how their outputs compare
if more than one of them works.

Examples~\ref{example:euler-n=3}, \ref{example:symmetry-n=3}, 
\ref{example:recip-n=3}, and \ref{example:K3-M(G)}
illustrate the methods in the $n=3$ case, and provide the first example 
in which all of them break down: when $a=12$, the denominator in the
formula for $h$ vanishes. 
(For $a>12$, we get $h=\frac{2a}{12-a}<0$, 
a negative ``fake Coxeter number.'') 
Here is another such example. 

\begin{example}
\label{example:K4} 

Let $G\!=\!K_4$ be the complete graph on $4$ vertices,
with all $6$ edges labeled by~$3$. 
Then each $3$-vertex subgraph $G'\lessdot G$
is of affine type~$\tilde A_2$, so \eqref{eq:M(G),n=3} gives
\hbox{$M(G')=3$}. 
Then \eqref{eq:M(G)-ingredients} yields $\Sigma_1=4\cdot 3=12$, 
and the denominator in
\eqref{eq:M(G)-rec} vanishes.
So the method based on Lemma~\ref{lem:M(G)-recur} does not work for
$G=K_4$. 
In fact, none of the methods described above will work in this case. 
Consequently, the methods will fail for 
any Coxeter diagram that has $K_4$ as an induced subgraph. 

Another similar example is the $4$-cycle with edge labels $3,4,3,4$
(in this order). 
\end{example}

Actually, each of our three main methods (see
Sections~\ref{sec:euler-char-method}--\ref{sec:reciprocity-based}) 
will typically fail in less drastic ways. 
In the Euler characteristic method and the (general) reciprocity-based method,
solving for $h$ can yield a non-constant rational function of~$m$;
see \eqref{eq:euler-char-phi-G-recur} and \eqref{eq:recip-method-7}, respectively.
In the symmetry-based method, the polynomial $Q(G,m)$ can 
fail to be symmetric with respect to the involution~\eqref{eq:mu-reflection}.  
In these cases, the methods yield an 
answer which should  be viewed with some skepticism.
In what follows, references to the failure of one of the three methods
mean that $h$ is non-constant or that symmetry fails.
(See Remark~\ref{rem:recip-special} for a discussion of the special 
versions of the reciprocity-based method, described in
Sections~\ref{sec:reciprocity-based}--\ref{sec:M(G)}.)

We now summarize the results of applying our methods to 
various Coxeter diagrams of infinite type, including the affine
diagrams.  The only infinite family for which the methods 
succeed wonderfully is the affine series of type~$\tilde A$.

\pagebreak[3]

\begin{proposition}
\label{prop:affineA-fake}
For a Coxeter diagram of type $\tilde A_{n-1}$ with $n\ge 3$, 
the fake invariants arising from
all methods described in Section~\ref{sec:coxeter-inv} 
agree, and coincide 
with the corresponding invariants of type~$B_{n}$.
That is, $h=2n$ and the exponents are $1,3,\ldots 2n-1$. 
\end{proposition} 

\begin{proof} 
Let the Coxeter diagram $G$ be a ``simply-laced $n$-cycle,'' 
i.e., an unoriented graph with the vertices $1,2,\dots,n$ such that the labels of
the edges $\{1,2\},\dots,\{n\!-\!1,n\}, \{n,1\}$ are equal to~$3$,
and the labels of all other edges are~$2$. 
To prove Proposition~\ref{prop:affineA-fake}, one needs to verify the
following statements for this particular diagram~$G$ and the invariants
given by the type-$B_n$ formulas: 
the recurrence~\eqref{eq:diffeq-f}; 
the formula \eqref{eq:euler-char-phi-G} for the ``fake Euler
characteristic;'' and the inclusion-exclusion
formula~\eqref{eq:recip-method-3}. 
This is a straightforward verification, omitted here. 
(Note that each connected induced subgraph of~$G$ is of type~$A$.)
\end{proof}

\begin{remark}
The recurrences that determine the fake invariants of type $\tilde
A_{n-1}$ are different from the type-$B_n$ recurrences,
so the coincidence of the final answers is somewhat mysterious. 
On the other hand, the fact that the $n$-cycle is combinatorially
related to the type-$B_n$ cluster complex (the $n$-dimensional
cyclohedron) is well known: there is a \emph{different} generalization
of associahedra \cite{carr-devadoss, DTS, postnikov-asso} 
which produces the cyclohedron in case of an $n$-cycle. 
It would be interesting to clarify the relationship between the two
constructions, and provide a conceptual explanation for
Proposition~\ref{prop:affineA-fake}. 
\end{remark}

Our computations suggest that the complete list of connected Coxeter
diagrams  for which all methods described in
Sections~\ref{sec:euler-char-method}--\ref{sec:reciprocity-based} 
produce \emph{positive integer} (fake) exponents is as follows:
\begin{itemize}
\item
Coxeter diagrams of finite irreducible Coxeter groups; 
\item
Coxeter diagrams of affine type $\tilde A$; 
\item
diagrams of rank~$3$ with the sum of edge labels equal to
$8$, $9$, $10$, or~$11$. 
\end{itemize} 

\begin{remark}
\label{rem:recip-special}
The discussion below does not refer to the methods based on
\eqref{eq:recip-method-1}--\eqref{eq:recip-method-3} 
and on Lemma~\ref{lem:M(G)-recur}, respectively. 
This is because these methods produce results which are, by construction,
the specializations (at $m=1$ and $m=0$, respectively)
of the results delivered by the general reciprocity method 
of Section~\ref{sec:reciprocity-based} 
based on the recurrences \eqref{eq:recip-method-5}
and~\eqref{eq:recip-method-7}. 
So whenever the general reciprocity method produces a fake Coxeter number 
$h$ that is a non-constant rational function in~$m$, 
the specializations of this rational function at $m=0$ and $m=1$ 
should be viewed with suspicion. 
Still, it is tempting to mention some instances of Coxeter diagrams of
affine type where the general
reciprocity method fails whereas the $M(G)$-based method of
Section~\ref{sec:M(G)} gives a positive integer value
of~$h$. 
They are: 
\begin{itemize}
\item 
type~$\tilde B_5\,$, with~$h=22$ and $M(G)=26$ (in type~$\tilde B_4\,$,
$h$ is a fraction);
\item
type~$\tilde C_n\,$, with $h=3n+4$  and $M(G)=3n+2$ (but irrational exponents); 
\item
type~$\tilde E_8\,$, with $h=98$ and $M(G)=306$. 
\end{itemize}
\end{remark}

The fake Coxeter number, if defined, is always a rational number (by
design), and sometimes even an integer. 
On the other hand, in most cases, some of the fake exponents are
irrational, 
as the generalized Fuss-Catalan invariant $N(G,m)$ would not
factor into linear polynomials in~$m$ (over~$\mathbb{Q}$).  
Still, $N(G,m)$ would sometimes be an integer-valued polynomial, 
potentially allowing for a meaningful interpretation. 

For the affine types other than~$\tilde A$, calculations give the
following results. 

\pagebreak[3]

The types $\tilde B_2/\tilde C_2$ and $\tilde{G}_2$ fall into the
$n=3$ case mentioned above, with $a=10$ and $a=11$, respectively.  
In type~$\tilde B_2/\tilde C_2\,$, we get $h=10$, with exponents
$1,5,9$---the same invariants as in type~$H_3$. 
In type~$\tilde{G}_2\,$, we get 
$h=22$, with exponents $1,11,21$.

In type $\tilde{B}_3$, all methods yield $h=\frac{76}{5}$, with exponents 
$1,\frac{33}{5},\frac{43}{5},\frac{71}{5}$.

In type $\tilde{C}_3$, all methods yield $h=13$, with exponents 
$1,\frac{13-\sqrt{17}}{2},\frac{13+\sqrt{17}}{2},12$.
Cf.\ Example~\ref{example:4-chain}
with $b_1=4$, $b_2=3$, $b_3=4$.

The type $\tilde{D}_4$ is the only ambiguous case that 
we have found.  
The reciprocity method gives $h=14$ with the exponents $1,6,6,9,13$.
This sequence of exponents is not symmetric, and
accordingly, the symmetry-based method fails.
However, the answer produced by blindly applying the formulas
from the symmetry-based  method
still agrees with the result of the reciprocity method.
The Euler characteristic method fails.

For all other affine types,  
computations up to rank 12 suggest that all methods~fail.

Beyond rank~$3$, our tests succeeded 
in a few isolated examples of non-affine infinite type, 
although they did not ordinarily produce integer values of~$h$.
Here are a few typical examples. 

For a $4$-cycle with one edge labeled~$4$ 
and the others labeled~$3$, 
we get $h=43/2$ with exponents 
$1,\frac{43-\sqrt{145}}{4},\frac{43+\sqrt{145}}{4},\frac{41}{2}$.
For a $5$-cycle with similar labels, all tests fail. 
For a $4$-cycle with one edge labeled~$5$ 
and the others labeled~$3$, all methods yield
 $h=-22$ with exponents 
$1,-11+2\sqrt{3},-11-2\sqrt{3},-23$.
We also found two more examples with positive integer fake Coxeter numbers: 
\[
\begin{array}{lll}
\bullet\overunder{3}{}\bullet\overunder{4}{}\bullet\overunder{4}{}\bullet
&\qquad & h=98\\
\bullet\overunder{4}{}\bullet\overunder{3}{}\bullet\overunder{5}{}\bullet
&\qquad & h=104
\end{array}
\]

We note that in all cases we have tried, 
if all methods work, then they all agree.

\pagebreak[3]

\section*{Acknowledgments}

We are grateful to Christos Athanasiadis for valuable comments on 
the preliminary versions of this paper, for explaining some of his
work, and for making available his paper~\cite{Ath-Tz} 
(joint with E.~Tzanaki) which we cite in Section~\ref{sec:recip}. 

Christian Krattenthaler, David Speyer and Hugh Thomas contributed 
several suggestions which improved the quality of the final version.
We also thank Egon Schulte, Richard Stanley, John
Stembridge, and Andrei Zelevinsky for helpful conversations. 

The first version of this paper was circulated at the workshop
``Braid groups, clusters and free probability'' held at the American
Institute of Mathematics on January 10--14, 2005.
We thank the organizers (Jon McCammond, Alexandru Nica, and Victor
Reiner) as well as the rest of participants for
stimulating discussions. 

This paper was completed in the spring of 2005, when S.F.\ was visiting the Mittag-Leffler
Institute in Djursholm, Sweden. 

Our computations of (fake) Coxeter
invariants were done with \texttt{Maple} 
using John Stembridge's package \texttt{coxeter}.
The \texttt{Maple} source code is currently available at 
\texttt{http://www.math.lsa.umich.edu/\textasciitilde nreading/papers/gcccc/}.

\end{document}